\documentclass[hidelinks,onefignum,onetabnum]{siamart220329}
\usepackage{amsmath}
\usepackage{amssymb}
\usepackage{enumerate}
\usepackage{mathrsfs}
\usepackage{cases}
\usepackage{color}
\usepackage{tikz}
\usepackage{mathtools}
\usepackage{mdframed}
\usepackage{makecell}
\usepackage{multirow}
\usepackage{siunitx}
\usepackage{graphicx}
\definecolor{darkblue}{rgb}{0.0, 0.0, 0.55}
\definecolor{bordeaux}{rgb}{0.34, 0.01, 0.1}
\newsiamthm{example}{Example}

\usepackage{lipsum}
\usepackage{amsfonts}
\usepackage{graphicx}
\usepackage{epstopdf}
\usepackage{algorithmic}
\ifpdf
  \DeclareGraphicsExtensions{.eps,.pdf,.png,.jpg}
\else
  \DeclareGraphicsExtensions{.eps}
\fi

\newsiamremark{remark}{Remark}
\newsiamremark{hypothesis}{Hypothesis}
\crefname{hypothesis}{Hypothesis}{Hypotheses}
\newsiamthm{claim}{Claim}
\newtheorem{conjecture}[theorem]{Conjecture}

\headers{Strengthened Complex Moment Hierarchies}{Jie Wang}

\title{Strengthened Complex Moment Hierarchies and Finite Convergence over Spheres
}

\author{Jie Wang\thanks{State Key Laboratory of Mathematical Sciences, Academy of Mathematics and Systems Science, Chinese Academy of Sciences, Beijing, China
  (\email{wangjie212@amss.ac.cn}, \url{https://wangjie212.github.io/jiewang/})}}

\usepackage{amsopn}

\def\R{{\mathbb{R}}}
\def\C{{\mathbb{C}}}
\def\N{{\mathbb{N}}}

\def\z{{\mathbf{z}}}

\def\y{{\mathbf{y}}}
\def\bu{{\mathbf{u}}}
\def\bw{{\mathbf{w}}}

\def\be{{\mathbf{e}}}
\def\bp{{\mathbf{p}}}

\def\a{{\boldsymbol{\alpha}}}

\def\b{{\boldsymbol{\beta}}}
\def\g{{\boldsymbol{\gamma}}}

\def\bv{{\boldsymbol{v}}}

\def\B{{\mathscr{B}}}

\def\M{{\mathbf{M}}}

\def\i{\hbox{\bf{i}}}

\def\rank{\hbox{\rm{rank}}}

\begin{document}

\maketitle

\begin{abstract}
This paper proposes strengthened complex moment relaxations for complex polynomial optimization. The proposed relaxation augments the pruned complex moment relaxation with selected positive semidefinite blocks from the full complex moment matrix. For the first strengthening, these blocks are motivated by the normality of multiplication operators and provide a tunable tradeoff between relaxation strength and computational cost through an additional normal order. In analogy with the real moment hierarchy, a flat-truncation condition for detecting global optimality is derived for the strengthened hierarchy. Furthermore, by considering the joint normality of multiplication operators, we propose a second strengthened complex moment hierarchy. For this stronger hierarchy, we are able to prove generic finite convergence and flat truncation for optimization over spheres. Numerical experiments demonstrate that the proposed strengthenings recover almost all the strength of the full relaxation at substantially lower computational cost.
\end{abstract}



\section{Introduction}
In this paper, we are concerned with the complex polynomial optimization problem
\begin{equation}\label{cpop}
f_{\min}\coloneqq\inf\,\{f(\z,\bar \z): \z\in \mathbf{K}\},\tag{CPOP}
\end{equation}
where
\begin{equation}
\mathbf{K}\coloneqq\{\z\in \mathbb C^n\mid g_i(\z,\bar \z)\geq 0,\ i=1,\ldots,m,\, h_j(\z,\bar \z)=0,\ j=1,\ldots,l\},
\end{equation}
and $f,g_1,\ldots,g_m,h_1,\ldots,h_l\in \mathbb C[\z,\bar \z]$ are self-conjugate polynomials, so that they take real values on $\mathbb C^n$. Such problems arise naturally in various fields, e.g., signal processing \cite{aubry2013ambiguity,dumitrescu2007positive,fogel2016phase,mariere2003,waldspurger2015phase}, power systems \cite{bienstock2020,josz2018lasserre}, quantum information \cite{fang2021sum,gribling2022bounding,le2025flat}, combinatorial optimization \cite{goemans2001approximation,sinjorgo2024cuts,zhang2006complex}.

A standard way to solve complex polynomial optimization is to convert each complex variable into its real and imaginary parts and then apply the real moment-SOS hierarchy of Lasserre \cite{Las01}. This approach is powerful and enjoys strong theoretical guarantees, but it doubles the number of variables and may become computationally expensive. To avoid this explosion, Josz and Molzahn \cite{josz2018lasserre} introduced a complex moment hierarchy that works directly with complex variables, where the moment matrix is indexed only by holomorphic monomials. This complex moment hierarchy has much lower computational complexity than the corresponding real moment hierarchy and has been shown to be effective in several applications, especially in optimal-power-flow problems \cite{josz2018lasserre,wang2022exploiting}.

However, the computational advantage of the complex moment hierarchy comes at a price. The pruned complex moment matrix keeps only a subset of the moment information contained in the full complex moment matrix that is indexed by general mixed monomials in \(\z\) and \(\bar \z\). Consequently, the complex moment hierarchy may converge more slowly than the real hierarchy and may fail to achieve global optimality at low relaxation orders \cite{wang2025real}. The full complex moment relaxation, which is equivalent to the real moment relaxation after the real-imaginary reformulation, is stronger but significantly more expensive. This motivates the following question:

\vspace{0.5em}
\begin{mdframed}
\textit{Can one strengthen the complex moment hierarchy by adding selected low-cost constraints from the full complex moment hierarchy?}
\end{mdframed}
\vspace{0.5em}

The purpose of this paper is to develop such a strengthening. The basic idea is to augment the usual complex moment relaxation by positive semidefinite (PSD) constraints of the form
\begin{equation}\label{coptcon}
\mathbf{N}^i_s(\y)\coloneqq\begin{bmatrix}
\M^{\C}_{s}(\y)&\M^{\C}_{s}(z_i\y)\\
\M^{\C}_{s}(\bar{z}_i\y)&\M^{\C}_{s}(|z_i|^2\y)\\
\end{bmatrix}\succeq0.
\end{equation}
Here \(s\) is a tunable parameter.
The matrices \(\mathbf{N}^i_s(\y)\) are principal submatrices of an appropriate full complex moment matrix, but they are not chosen arbitrarily. They arise naturally from the normality of multiplication operators associated with multiplication by the complex variables \(z_i\). This operator-theoretic interpretation is central. In the real moment hierarchy, multiplication operators are symmetric and hence normal.
A single flat-truncation condition is sufficient for certifying global optimality and minimizer extraction. In the pruned complex setting, the analogous multiplication operators need not be normal. The PSD blocks above provide a controlled way to impose this missing normality information, while avoiding the heavy cost of the real or full complex hierarchy.


The proposed hierarchy therefore lies between the pruned complex hierarchy and the full complex, equivalently real, moment hierarchy. It is stronger than the pruned hierarchy because it imposes additional PSD constraints. It is much cheaper than the full complex moment hierarchy because it uses only selected small blocks of the full moment matrix. The parameter \(s\) provides an additional degree of flexibility: small values of \(s\) often improve the lower bound substantially with little additional cost, while larger values of \(s\) give stronger relaxations at the price of larger PSD blocks.
\vspace{0.5em}

The main contributions of this paper are summarized as follows.

First, we introduce a strengthened complex moment relaxation obtained by augmenting the pruned complex relaxation with normality-inspired PSD constraints. This gives a bilevel hierarchy indexed by the relaxation order $r$ and an additional parameter $s$. We also give an explicit example showing that the strengthened relaxation can be strictly tighter than the pruned complex relaxation at the same relaxation order. 

Second, we explain the structural reason for the added constraints. We show that, under the usual flatness assumptions, the conditions \(\forall i,\,\mathbf{N}^i_s(\y)\succeq0\) characterize the normality of multiplication operators. This gives a principled justification for choosing these particular submatrices of the full complex moment matrix rather than arbitrary submatrices.

Third, we show that the strengthened hierarchy admits a flat-truncation certificate for global optimality. Once the flatness condition is detected, one can extract global minimizers from the moment matrix, just as in the real moment hierarchy. This provides a simple stopping criterion that is unavailable for the pruned complex hierarchy alone.

Forth, we prove a generic finite convergence result with a stronger strengthening. The individual normality blocks \(\forall i,\,\mathbf{N}^i_s(\y)\succeq0\) enforce the normality of multiplication operators independently. We also study a larger PSD condition encoding the joint normality of multiplication operators. This leads to a second strengthened complex moment hierarchy. For optimization of a self-conjugate polynomial over the unit sphere, we prove that this joint-normality-inspired hierarchy has finite convergence and flat truncation under the generic assumptions that the set of global minimizers is nonempty and finite and that every global minimizer is nondegenerate on the sphere. The proof uses a special certificate for nonnegative Hermitian polynomials on the sphere, together with a matrix-valued holomorphic factorization argument.

Finally, we present numerical experiments on several classes of complex polynomial optimization problems. These include random quartic problems with unit-modulus variables, random quartic problems on spheres, sparse quartic problems on multi-spheres, and complex polynomial optimization formulations related to Smale's mean value conjecture and the Mordell inequality conjecture. The experiments show that the strengthened hierarchies can improve bounds without increasing the main relaxation order and can often certify global optimality with substantially less computational effort than higher-order pruned complex relaxations or full complex relaxations.

The paper is organized as follows. Section~2 introduces notation and recalls the pruned complex and full complex moment relaxations. Section~3 presents the strengthened complex moment hierarchy based on individual normality conditions, gives strict-improvement examples, explains the link with the normality of multiplication operators, and describes the sparse variant. Section~4 introduces the joint-normality-inspired strengthening and proves generic finite convergence and flat truncation for optimization over spheres. Section~5 reports numerical experiments. Section~6 concludes with an open question on finite convergence for more general constrained complex polynomial optimization problems.

\section{Notation and preliminaries}
Let $\N$ (resp.\,$\N^*$) be the set of nonnegative (resp.\,positive) integers.
For $n\in\N^*$, let $[n]\coloneqq\{1,2,\ldots,n\}$. For $\a=(\alpha_i)\in\N^n$, let $|\a|\coloneqq\sum_{i=1}^n\alpha_i$. For $r\in\N$, let $\N^n_r\coloneqq\{\a\in\N^n\mid|\a|\le r\}$. We use $A\succeq0$ to indicate that the matrix $A$ is PSD.
Let $\bar{a}$ denote the conjugate of a complex number $a$ and $\bv^{*}$ (resp. $A^{*}$) denote the conjugate transpose of a complex vector $\bv$ (resp. a complex matrix $A$). Let $\z=(z_1,\ldots,z_n)$ be a tuple of complex variables and $\bar{\z}=(\bar{z}_1,\ldots,\bar{z}_n)$ be the conjugate of $\z$. 
We denote by $\C[\z]\coloneqq\C[z_1,\ldots,z_n]$ (resp. $\C[\z,\bar{\z}]\coloneqq\C[z_1,\ldots,z_n,\bar{z}_1,\ldots,\bar{z}_n]$) the complex polynomial ring in $\z$ (resp. $\z,\bar{\z}$). Let $\C[\z]_r\coloneqq\left\{p=\sum_{\b}p_{\b}\z^{\b}\in\C[\z]\,\middle|\,|\b|\le r\right\}$. A polynomial $f\in\C[\z,\bar{\z}]$ can be written as $f=\sum_{(\b,\g)\in\N^n\times\N^n}f_{\b,\g}\z^{\b}\bar{\z}^{\g}$ with $f_{\b,\g}\in\C$. 
The \emph{conjugate} of $f$ is defined as $\bar{f}=\sum_{(\b,\g)\in\N^n\times\N^n}\bar{f}_{\b,\g}\z^{\g}\bar{\z}^{\b}$. The polynomial $f$ is \emph{self-conjugate} if $\bar{f}=f$. It is clear that self-conjugate polynomials take only real values. For $d\in\N$, $[\z]_{d}$ (resp.\,$[\z,\bar{\z}]_{d}$) stands for the standard monomial basis in $\z$ (resp.\,$\z,\bar{\z}$) of degree up to $d$. Let $\delta_{\bv}$ denote the Dirac measure centered at a point $\bv\in\C^n$. We use $\|\cdot\|$ to denote the $2$-norm of a vector.

By invoking Borel measures, \eqref{cpop} admits the following equivalent reformulation:
\begin{equation}\label{meas}
\begin{cases}
\inf\limits_{\mu\in\mathcal{M}_+(\mathbf{K})} &\int_{\mathbf{K}} f\,\mathrm{d}\mu\\
\quad\,\,\rm{s.t.}&\int_{\mathbf{K}}\,\mathrm{d}\mu=1,\\
\end{cases}
\end{equation}
where $\mathcal{M}_+(\mathbf{K})$ denotes the set of finite positive Borel measures on $\mathbf{K}$. 

Suppose that $\y=(y_{\b,\g})\in\C^{\N^n\times\N^n}$ is a complex sequence satisfying $y_{\b,\g}=\bar{y}_{\g,\b}$.\footnote{We always assume that a complex sequence $\y$ has this property in this paper.} We associate it with a linear functional $L_{\y}:\C[\z,\bar{\z}]\rightarrow\C$ by
\begin{equation*}
p=\sum_{(\b,\g)}p_{\b,\g}\z^{\b}\bar{\z}^{\g}\longmapsto L_{\y}(p)=\sum_{(\b,\g)}p_{\b,\g}y_{\b,\g}.
\end{equation*}
For any $r\in\N$, the $r$-th order \emph{(pruned) complex moment matrix} $\M^{\C}_{r}(\y)$ is the Hermitian matrix defined by
\begin{equation*}
\M^{\C}_r(\y)=L_{\y}\left([\z]_r\cdot[\z]_r^*\right)
\end{equation*}
and the $r$-th order \emph{full complex moment matrix} $\M^{2\C}_{r}(\y)$ is the Hermitian matrix defined by
\begin{gather*}
\M^{2\C}_r(\y)=L_{\y}\left([\z,\bar{\z}]_r\cdot[\z,\bar{\z}]_r^*\right).
\end{gather*}
Moreover, for a self-conjugate polynomial $g$,
the $r$-th order \emph{(pruned) complex localizing matrix} $\M^{\C}_{r}(g\y)$ associated with $g$ is the Hermitian matrix defined by
\begin{equation*}
\M^{\C}_r(g\y)=L_{\y}\left([\z]_r\cdot[\z]_r^*\cdot g\right)
\end{equation*}
and the $r$-th order \emph{full complex localizing matrix} $\M^{2\C}_{r}(g\y)$ associated with $g$ is the Hermitian matrix defined by
\begin{gather*}
\M^{2\C}_{r}(g\y)=L_{\y}\left([\z,\bar{\z}]_r\cdot[\z,\bar{\z}]_r^*\cdot g\right).
\end{gather*}


Assume $f=\sum_{(\b,\g)}f_{\b,\g}\z^{\b}\bar{\z}^{\g}$, $g_i=\sum_{(\b,\g)}g^i_{\b,\g}\z^{\b}\bar{\z}^{\g}$, and $h_j=\sum_{(\b,\g)}h^j_{\b,\g}\z^{\b}\bar{\z}^{\g}$. Let $d^f\coloneqq\max\,\{|\b|,|\g|: f_{\b,\g}\ne0\}$, $d^g_i\coloneqq\max\,\{|\b|,|\g|: g^i_{\b,\g}\ne0\}$ for $i\in[m]$, $d^h_j\coloneqq\max\,\{|\b|,|\g|: h^j_{\b,\g}\ne0\}$ for $j\in[l]$, and $d_{\min}\coloneqq\max\,\{d^f,d^g_1,\ldots,d^g_m,d^h_1,\ldots,d^h_l\}$. For a positive integer $r\ge d_{\min}$, the $r$-th order (pruned) complex moment relaxation \cite{josz2018lasserre} for \eqref{cpop} is given by
\begin{equation}\label{cmom}
\tau_r\coloneqq\begin{cases}
\inf\limits_{\y}& L_{\y}(f)\\
\rm{s.t.}&\M^{\C}_{r}(\y)\succeq0,\quad y_{\mathbf{0},\mathbf{0}}=1,\\
&\M^{\C}_{r-d^g_i}(g_i\y)\succeq0,\quad i\in[m],\\
&\M^{\C}_{r-d^h_j}(h_j\y)=0,\quad j\in[l].
\end{cases}
\end{equation}

Let $\tilde{d}^f\coloneqq\max\,\{|\b|+|\g|: f_{\b,\g}\ne0\}$, $\tilde{d}^g_i\coloneqq\max\,\{|\b|+|\g|: g^i_{\b,\g}\ne0\}$ for $i\in[m]$, $\tilde{d}^h_j\coloneqq\max\,\{|\b|+|\g|: h^j_{\b,\g}\ne0\}$ for $j\in[l]$, and \[\tilde{d}_{\min}\coloneqq\max\,\{\lceil\tilde{d}^f/2\rceil,\lceil\tilde{d}^g_1/2\rceil,\ldots,\lceil\tilde{d}^g_m/2\rceil,\lceil\tilde{d}^h_1/2\rceil,\ldots,\lceil\tilde{d}^h_l/2\rceil\}.\] For a positive integer $r\ge\tilde{d}_{\min}$, the $r$-th order full complex moment relaxation for \eqref{cpop} is given by
\begin{equation}\label{fcmom}
\rho_r\coloneqq\begin{cases}
\inf\limits_{\y}& L_{\y}(f)\\
\rm{s.t.}&\M^{2\C}_{r}(\y)\succeq0,\quad y_{\mathbf{0},\mathbf{0}}=1,\\
&\M^{2\C}_{r-\lceil\tilde{d}^g_i/2\rceil}(g_i\y)\succeq0,\quad i\in[m],\\
&L_\y(\z^{\b}\bar\z^{\g}h_j)=0,\quad \forall|\b|+|\g|\le 2r-\tilde{d}_j^h, \quad j\in[l].
\end{cases}
\end{equation}

\begin{remark}
By introducing real variables for the real and imaginary parts of each complex variable, respectively, \eqref{cpop} could be converted into a real polynomial optimization problem. In fact, the $r$-th order full complex moment relaxation for \eqref{cpop} is equivalent to the $r$-th order real moment relaxation for the equivalent real polynomial optimization problem; see \cite{lasserre2008semidefinite}.
\end{remark}

\begin{remark}
We refer the reader to \cite{wang2026more} for reformulating the complex semidefinite program (SDP) \eqref{s-cmom} as a real SDP.
\end{remark}

We say that the sequence $\y$ admits a (finitely atomic) representing measure if it can be realized by a Borel (finitely atomic) measure $\mu$, i.e., $y_{\b,\g}=\int_{\mathbf{K}} \z^{\b}\bar{\z}^{\g}\,\mathrm{d}\mu$ for any $\b,\g\in\N^n$. Accordingly, $\y$ is called the moment sequence of $\mu$. It is clear that the relaxations \eqref{cmom} and \eqref{fcmom} attain global optimality for \eqref{cpop} if the sequence $\y$ admits a representing measure.



\section{Strengthening by the normality conditions}\label{sec3}
For each $i\in[n]$, let $b^i_s(\z)$ be the column vector formed by concatenating $[\z]_s$ and $\bar{z}_i[\z]_s$. We have 
\begin{equation}
\mathbf{N}^i_s(\y)=\begin{bmatrix}
\M^{\C}_{s}(\y)&\M^{\C}_{s}(z_i\y)\\
\M^{\C}_{s}(\bar{z}_i\y)&\M^{\C}_{s}(|z_i|^2\y)\\
\end{bmatrix}=L_{\y}\left(b^i_s(\z)b^i_s(\z)^*\right).
\end{equation}
As $b^i_s(\z)$ is a subvector of $[\z,\bar{\z}]_{s+1}$,
$\mathbf{N}^i_s(\y)$ is a principal submatrix of $\M_{s+1}^{2\C}(\y)$.
We then propose to strengthen the pruned complex moment relaxation \eqref{cmom} with the PSD constraints \eqref{coptcon}:
\begin{equation}\label{s-cmom}
\tau'_{r,s}\coloneqq\begin{cases}
\inf\limits_{\y}& L_{\y}(f)\\
\rm{s.t.}&\M^{\C}_{r}(\y)\succeq0,\quad y_{\mathbf{0},\mathbf{0}}=1,\\
&\M^{\C}_{r-d^g_i}(g_i\y)\succeq0,\quad i\in[m],\\
&\M^{\C}_{r-d^h_j}(h_j\y)=0,\quad j\in[l],\\
&\mathbf{N}^i_s(\y)\succeq0,\quad i\in[n],
\end{cases}
\end{equation}
where $s\in\N^*$ is a tunable parameter.

\begin{proposition}\label{sec5:thm2}
The following are true for any $r\ge d_{\min}$ and $s\in\N^*$:
\begin{enumerate}[\rm (i)]
    \item $\tau_r\le\tau'_{r,s}\le\tau'_{r,s+1}\le f_{\min}$;
    \item $\tau'_{r,s}\le\tau'_{r+1,s}$;
    \item $\tau'_{r,s}\le\rho_{\max\,\{r,s+1\}}$.
\end{enumerate}
\end{proposition}
\begin{proof}
It is straightforward from the constructions.
\end{proof}


\begin{remark}
By Theorem 3.2 of \cite{josz2018lasserre}, the complex moment hierarchy \eqref{cmom} has asymptotic convergence (i.e., $\lim_{r\to\infty}\tau_r=f_{\min}$) when a sphere constraint is present. As an immediate corollary, the strengthened complex moment hierarchy \eqref{s-cmom} also has asymptotic convergence (i.e., $\lim_{r\to\infty}\tau'_{r,s}=f_{\min}$) under the presence of a sphere constraint.
\end{remark}

\begin{remark}
The conditions \eqref{coptcon} have been used in \cite[Corollary 4.1]{wang2025real} for detecting global optimality of the relaxation \eqref{cmom}.
Similar but more complicated conditions have appeared in \cite[Proposition 4.1]{josz2018lasserre}.
\end{remark}

We give a simple example showing that the strengthened complex moment relaxation can be strictly tighter than the pruned complex moment relaxation at the same relaxation order.
\begin{example}
Consider the complex polynomial optimization problem
\begin{equation}
f_{\min}\coloneqq\inf_{\z\in\mathbb C^2}
\left\{f(\z,\bar \z)=|\bar z_1-z_2|^2:
|z_1|^2+|z_2|^2=1\right\}.
\end{equation}
Clearly, $f(\z,\bar \z)\geq 0$ on the feasible set. Moreover, $f_{\min}=0$ is attained, for instance, at $z_1=z_2=\frac{1}{\sqrt 2}$.
We now compare the second-order pruned complex moment relaxation with the strengthened relaxation of normal order $s=1$. Let $\{\be_1,\be_2\}$ be the standard basis vector of $\R^2$.
For the pruned complex moment relaxation at order $r=2$, the sphere constraint gives $\M_1^{\mathbb C}\big((1-|z_1|^2-|z_2|^2)\y\big)=0$, in particular,
\begin{equation*}
y_{\be_1,\be_1}+y_{\be_2,\be_2}=1,\,
y_{\be_1,\be_1}=y_{2\be_1,2\be_1}+y_{\be_1+\be_2,\be_1+\be_2},\,y_{\be_2,\be_2}=y_{\be_1+\be_2,\be_1+\be_2}+
y_{2\be_2,2\be_2}.
\end{equation*}
Since $\M_2^{\mathbb C}(\y)\succeq 0$, all diagonal entries are nonnegative. Therefore,
\begin{equation}
2y_{\be_1+\be_2,\be_1+\be_2}
\leq
y_{\be_1,\be_1}+y_{\be_2,\be_2}=1,
\end{equation}
and hence
$y_{\be_1+\be_2,\be_1+\be_2}\leq \frac12$.
Again by $\M_2^{\mathbb C}(\y)\succeq 0$, the Cauchy-Schwarz inequality for the entries indexed by $1$ and $z_1z_2$ gives
\begin{equation}
|y_{\be_1+\be_2,\mathbf{0}}|^2
\leq
y_{\be_1+\be_2,\be_1+\be_2}y_{\mathbf{0},\mathbf{0}}
\leq
\frac12.
\end{equation}
Thus, $\operatorname{Re}(y_{\be_1+\be_2,\mathbf{0}})
\leq\frac{1}{\sqrt 2}$.
The objective value of the pruned relaxation is
\begin{equation}
L_\y(f)=y_{\be_1,\be_1}+y_{\be_2,\be_2}-y_{\be_1+\be_2,\mathbf{0}}-
y_{\mathbf{0},\be_1+\be_2}=
1-2\operatorname{Re}(y_{\be_1+\be_2,\mathbf{0}}).
\end{equation}
Consequently, $L_\y(f)\geq 1-\sqrt 2$.
This lower bound is attained. Indeed, 
define
\begin{equation}
\M_2^{\mathbb C}(\y)
=
\begin{bmatrix}
1 & 0 & 0 & 0 & 1/\sqrt 2 & 0\\
0 & 1/2 & 0 & 0 & 0 & 0\\
0 & 0 & 1/2 & 0 & 0 & 0\\
0 & 0 & 0 & 0 & 0 & 0\\
1/\sqrt 2 & 0 & 0 & 0 & 1/2 & 0\\
0 & 0 & 0 & 0 & 0 & 0
\end{bmatrix}.
\end{equation}
This matrix is PSD, because its only nontrivial $2\times 2$ block is
\begin{equation}
\begin{bmatrix}
1 & 1/\sqrt 2\\
1/\sqrt 2 & 1/2
\end{bmatrix}\succeq 0.
\end{equation}
A direct verification also gives
$\M_1^{\mathbb C}\big((1-|z_1|^2-|z_2|^2)\y\big)=0$.
Therefore, this $\y$ is feasible for the second-order pruned complex moment relaxation and yields
$L_\y(f)=1-\sqrt 2$.
Hence $\tau_2=1-\sqrt 2<0$.

We now show that the strengthened relaxation with $r=2$ and $s=1$ is exact. The PSDness of $\mathbf{N}^1_1(\y)$ implies
$L_\y(f)=L_\y\big(|\bar z_1-z_2|^2\big)\geq 0$.
Therefore, every feasible point of the strengthened relaxation satisfies
$L_\y(f)\geq 0$.
On the other hand, the moment sequence of the feasible point $z_1=z_2=\frac{1}{\sqrt 2}$ is feasible for the strengthened relaxation and gives an objective value $0$. Hence $\tau'_{2,1}=0$.
Thus, we obtain
$\tau_2=1-\sqrt 2<0=\tau'_{2,1}=f_{\min}$.
\end{example}

\begin{example}
Let us revisit Example 4.1 of \cite{josz2018lasserre}:
\begin{equation}
\begin{cases}
\inf\limits_{z_1,z_2\in\C} &3-|z_1|^2-0.5\i z_1\bar{z}_2^2+0.5\i z_2^2\bar{z}_1\\
\,\,\,\,\,\mathrm{s.t.}&z_2+\bar{z}_2\ge0,\\
&|z_1|^2-0.25z_1^2-0.25\bar{z}_1^2-1=0,\\
&|z_1|^2+|z_2|^2-3=0,\\
&\i z_2-\i\bar{z}_2=0.
\end{cases}
\end{equation}
The second-order pruned complex moment relaxation yields a lower bound $0.155089$ and the third-order pruned relaxation yields the global optimum $0.428175$. In \cite{josz2018lasserre}, it was reported that by adding a $9\times9$ PSD constraint
\begin{equation}
\begin{bmatrix}
\M^{\C}_{1}(\y)&\M^{\C}_{1}(z_1\y)&\M^{\C}_{1}(z_2\y)\\
\M^{\C}_{1}(\bar{z}_1\y)&\M^{\C}_{1}(|z_1|^2\y)&\M^{\C}_{1}(z_2\bar{z}_1\y)\\
\M^{\C}_{1}(\bar{z}_2\y)&\M^{\C}_{1}(z_1\bar{z}_2\y)&\M^{\C}_{1}(|z_2|^2\y)
\end{bmatrix}\succeq0
\end{equation}
to the second-order pruned relaxation, one obtains the global optimum $0.428175$. Instead, by adding two $6\times6$ PSD constraints
\begin{equation}
\mathbf{N}^1_1(\y)\succeq0,\quad
\mathbf{N}^2_1(\y)\succeq0
\end{equation}
to the second-order pruned relaxation, one also obtains the global optimum $0.428175$.
\end{example}

In the following proposition, we show that if each complex variable has unit norm, then the strengthening \eqref{s-cmom} is effective only if $s\ge r$.
\begin{proposition}\label{unitnorm}
Suppose that \eqref{cpop} contains the unit-modulus constraint for each complex variable, then $\tau'_{r,s}=\tau_{r}$ for $s<r$.
\end{proposition}
\begin{proof}
The PSD constraint $\mathbf{N}^i_s(\y)\succeq0$ in \eqref{s-cmom} corresponds to the term $|p(\z)+\bar z_i q(\z)|^2$ with $p(\z),q(\z)\in\C[\z]_s$ in the dual problem. Since each $|z_i|=1$, we have $|p(\z)+\bar z_i q(\z)|^2=|z_i(p(\z)+\bar z_i q(\z))|^2=|z_ip(\z)+q(\z)|^2$. Because $z_ip(\z)+q(\z)\in\C[\z]_{s+1}$, the term $|z_ip(\z)+q(\z)|^2$ corresponds to a PSD block of $\M_{s+1}^{\C}(\y)\succeq0$ after taking the dual. Therefore, the constraint $\mathbf{N}^i_s(\y)\succeq0$ is effective in \eqref{s-cmom} only if $r\le s$.
\end{proof}

\subsection{Link to the normality of multiplication operators}
As there are numerous ways to pick constraints from the full complex moment relaxations,
one may wonder why we choose the particular conditions \eqref{coptcon} to strengthen the complex moment relaxations.
In the following, we show that those conditions arise naturally in ensuring the normality of multiplication operators and enable the strengthened complex moment hierarchy to share the same global optimality criterion with the real moment hierarchy.

For $p\in\C[\z]_{r}$, we write $\bp$ for the coefficient vector of $p$ such that $p=\bp^{\intercal}[\z]_{r}$. 
Let $\y$ be a complex sequence satisfying $\M^{\C}_{s}(\y)\succeq0$ and $\rank\,\M_{s}^{\C}(\y)=\rank\,\M^{\C}_{s+1}(\y)$. Define a sesquilinear form associated with $\y$ by $\langle p,q\rangle_\y=L_\y(p\bar q)$. Let $\ker\M_{s}^{\C}(\y)\coloneqq\left\{p\in\C[\z]_s\mid\M_{s}^{\C}(\y)\bp=\mathbf{0}\right\}$ and similarly for $\ker\M_{s+1}^{\C}(\y)$. One can easily show that \[\C[\z]_{s}/\ker\M_{s}^{\C}(\y)\cong\C[\z]_{s+1}/\ker\M_{s+1}^{\C}(\y).\]
The multiplication operators $M_i,i\in[n]$ are defined by
\begin{align}\label{multi-oper}
  M_i:  \C[\z]_{s}/\ker\M_{s}^{\C}(\y)&\longrightarrow\C[\z]_{s+1}/\ker\M_{s+1}^{\C}(\y), \\
  p\quad&\longmapsto\quad z_ip.\notag
\end{align}

We say that a linear operator $T$ is \emph{normal} if $[T^{*},T]\coloneqq T^{*}T-TT^{*}=0$ where $T^*$ denotes the adjoint operator of $T$. In case $T$ is of finite dimension, the normality of $T$ is equivalent to the hyponormality, i.e., $[T^{*},T]\succeq0$.

\begin{proposition}\label{normal}
Suppose that $\y$ is a complex sequence satisfying $\M^{\C}_{s}(\y)\succeq0$ and $\rank\,\M_{s}^{\C}(\y)=\rank\,\M^{\C}_{s+1}(\y)$.
Then for any $i\in[n]$, the multiplication operator $M_i$ defined in \eqref{multi-oper} is normal if and only if $\mathbf{N}^i_s(\y)\succeq0$.
\end{proposition}
\begin{proof}
Let $\mathcal{P}_s$ denote the orthogonal projection onto $\C[\z]_{s}/\ker\M_{s}^{\C}(\y)$. The adjoints of the multiplication operators satisfy $M_i^*(p)=\mathcal{P}_s(\bar zp), \,i\in[n]$. Fix any $q\in\C[\z]_s$ and let $\mathbf{q}$ be the coefficient vector of $q$ with respect to the monomial basis $z_i[\z]_s$.
Let us write $\mathbf{N}^i_s(\y)$ in block form
\begin{equation*}
    \mathbf{N}^i_s(\y)=\begin{bmatrix}
        A&B\\B^*&C
    \end{bmatrix}.
\end{equation*}
Assume that $A\succ0$. Otherwise, use the Moore–Penrose inverse $A^{\dagger}$ instead of $A^{-1}$ in the following. The Schur complement of $A$ is $C-B^*A^{-1}B$, and we have
\begin{align*}
    \mathbf{q}^*(C-B^*A^{-1}B)\mathbf{q}
    &=\langle \bar z_iq,\bar z_iq\rangle_\y-\left\langle \mathcal{P}_s(\bar z_iq),\mathcal{P}_s(\bar z_iq)\right\rangle_\y\\
    &=\langle z_iq,z_iq\rangle_\y-\langle M_i^*(q),M_i^*(q)\rangle_\y\\
    &=\langle M_i(q), M_i(q)\rangle_\y-\langle q,M_iM_i^*(q)\rangle_\y\\
    &=\langle q, M_i^*M_i(q)\rangle_\y-\langle q,M_iM_i^*(q)\rangle_\y\\
    &=\langle q,[M_i^*,M_i](q)\rangle_\y.
\end{align*}
We thus conclude that $\mathbf{N}^i_s(\y)\succeq0$ if and only if $\M^{\C}_{s}(\y)\succeq0$ and $[M_i^*,M_i]\succeq0$, which implies the desired conclusion.
\end{proof}

In view of Proposition \ref{normal}, we call the index $s$ in \eqref{s-cmom} the \emph{normal order}. Finite convergence of the strengthened complex moment hierarchy could be detected in the same way as for the real moment hierarchy \cite{nie2013certifying}, which is stated in the following theorem.
Let $d_{\mathbf{K}}\coloneqq\max\,\{d^g_1,\ldots,d^g_m,d^h_1,\ldots,d^h_l\}$.
\begin{theorem}[Global optimality]\label{opt}
Suppose that $\y$ is an optimal solution of the strengthened relaxation \eqref{s-cmom}. If there is an integer $t$ such that $d_{\min}\leq t \leq\min\,\{s+d_{\mathbf{K}}, r\}$ and $\rank\,\M^{\C}_t(\y) = \rank\,\M^{\C}_{t-d_{\mathbf{K}}}(\y)$, then the relaxation \eqref{s-cmom} for \eqref{cpop} is tight, i.e., $\tau'_{r,s} = f_{\min}$.
\end{theorem}
\begin{proof}
The flatness condition allows us to define the multiplication operators $M_i:\C[\z]_{t-d_{\mathbf{K}}}/\ker\M_{t-d_{\mathbf{K}}}^{\C}(\y)\rightarrow\C[\z]_{t-d_{\mathbf{K}}+1}/\ker\M_{t-d_{\mathbf{K}}+1}^{\C}(\y),\,i\in[n]$ which are  pairwise commuting. By Proposition \ref{normal}, each $M_i$ is normal. Hence the $M_i$'s are simultaneously unitarily diagonalizable. So there exist $a\coloneqq\rank\,\M^{\C}_t(\y)$ points $\{\bu^{(1)},\cdots,\bu^{(a)}\}$
and an orthonormal basis $q_1,\ldots,q_a$ of $\C[\z]_{t-d_{\mathbf{K}}}/\ker\M_{t-d_{\mathbf{K}}}^{\C}(\y)$ such that $M_i(q_j)=u_i^{(j)}q_j$.
Define weights $w_j=|\langle q_j,1\rangle_\y|^2\ge0$, $j\in[a]$. The equality $y_{\mathbf{0},\mathbf{0}}=1$ yields $\sum_{j=1}^aw_j=1$. For every pair of monomials $\z^\a,\z^\b$ of degree within the truncation, we have
\[L_{\y}(\z^\a\bar\z^\b)=\langle M^{\a}(1),M^{\b}(1)\rangle_\y=\sum_{j=1}^aw_j(\bu^{(j)})^\a(\overline{\bu^{(j)}})^\b,\] where $M^{\g}\coloneqq \prod_{i=1}^nM_i^{\gamma_i}$.
Thus, $\y$ has the finitely atomic representing measure $\mu=\sum_{j=1}^aw_j\delta_{\bu^{(j)}}$. Moreover, the constraints of localizing matrices enforce the representing measure $\mu$ is supported on the feasible set $\mathbf{K}$. Therefore, the relaxation is tight.
\end{proof}

\begin{remark}
The rank condition in Theorem \ref{opt} is called \emph{flat truncation}. Once it is satisfied, one could extract $\rank\,\M^{\C}_t(\y)$ global minimizers for \eqref{cpop} from the moment matrix \cite{klep2018minimizer}.
\end{remark}

\subsection{Integrating with correlative sparsity}
The strengthening approach can be further integrated with correlative sparsity to improve scalability. Suppose that the index sets $[n]$, $[m]$, $[l]$ can be decomposed into $\{I_1,\dots,I_p\}$, $\{J_1,\dots,J_p\}$, and $\{K_1,\dots,K_p\}$, respectively, such that 1) $f = f_1 + \cdots + f_p$ with $f_k \in\C[\z_{I_k},\bar{\z}_{I_k}]$ for $k\in[p]$; 2) for all $k\in[p]$ and $i\in J_k$, $g_i\in\C[\z_{I_k},\bar{\z}_{I_k}]$; 3) for all $k\in[p]$ and $j\in K_k$, $h_j\in\C[\z_{I_k},\bar{\z}_{I_k}]$, where $\C[\z_{I_k},\bar{\z}_{I_k}]$ denotes the polynomial ring in those variables indexed by $I_k$. Let $\M^{\C}_r(\y, I_k)$ (resp. $\M^{\C}_r(g\y, I_k)$) be the submatrix obtained from $\M^{\C}_r(\y)$ (resp. $\M^{\C}_r(g\y)$) by retaining only those rows and columns indexed by $\b=(\beta_i)\in\N_r^n$ of $\M^{\C}_r(\y)$ (resp. $\M^{\C}_r(g\y)$) with $\beta_i = 0$ if $i \notin I_k$. Then, we can strengthen the sparse complex moment relaxation for complex polynomial optimization as follows:
\begin{equation}
\begin{cases}
\inf\limits_{\y} &L_{\y}(f)\\
\rm{s.t.}&\M^{\C}_{r}(\y, I_k)\succeq0,\quad k\in[p],\\
&\M^{\C}_{r-d^g_i}(g_i\y, I_k)\succeq0,\quad i\in J_k, k\in[p],\\
&\M^{\C}_{r-d^h_j}(h_j\y, I_k)=0,\quad j\in K_k, k\in[p],\\
&\begin{bmatrix}
\M^{\C}_{s}(\y, I_k)&\M^{\C}_{s}(z_i\y, I_k)\\
\M^{\C}_{s}(\bar{z}_i\y, I_k)&\M^{\C}_{s}(|z_i|^2\y, I_k)\\
\end{bmatrix}\succeq0,\quad i\in I_k,k\in[p],\\
&y_{\mathbf{0},\mathbf{0}}=1.
\end{cases}
\end{equation}

\begin{remark}
One may also take term sparsity into account for the strengthened approach. We refer the reader to \cite{tssos1,tssos3} and omit the details for conciseness.
\end{remark}

\section{The joint normality condition and finite convergence over spheres}
Nie proved in his seminal work \cite{nie2013certifying,nie2014optimality,nie2023moment} that the real moment hierarchy has finite convergence and flat truncation in the generic case under the Archimedean condition. A natural question is: Can we establish similar results for the strengthened complex moment hierarchy? Although the normality conditions introduced in Section~\ref{sec3} are not strong enough to guarantee generic finite convergence, we shall prove in this section that, if the \emph{joint normality condition} of multiplication operators is used to strengthen the complex moment hierarchy, then both finite convergence and flat truncation generically occur for complex polynomial optimization over spheres.


We say that a tuple of operators $\{T_1,\ldots,T_n\}$ is \emph{jointly normal} if $[T_i,T_j]=0$ and $[T_i^*,T_j]=0$ for all $1\le i,j\le n$.


\begin{proposition}\label{joint-normal}
Suppose that $\y$ is a complex sequence satisfying $\M^{\C}_{s}(\y)\succeq0$ and $\rank\,\M_{s}^{\C}(\y)=\rank\,\M^{\C}_{s+1}(\y)$.
Then the multiplication operators $M_1,\ldots,M_n$ defined in \eqref{multi-oper} jointly normal if and only if
\begin{equation}
\mathbf{H}_{s}(\y)\coloneqq\begin{bmatrix}
\M^{\C}_{s}(\y)&\M^{\C}_{s}(\bar z_1\y)&\cdots &\M^{\C}_{s}(\bar z_n\y)\\
\M^{\C}_{s}(z_1\y)&\M^{\C}_{s}(|z_1|^2\y)&\cdots &\M^{\C}_{s}(z_1\bar z_n\y)\\
\vdots&\vdots &\ddots&\vdots\\
\M^{\C}_{s}(z_n\y)&\M^{\C}_{s}(z_n\bar z_1\y)&\cdots &\M^{\C}_{s}(|z_n|^2\y)
\end{bmatrix}\succeq0.
\end{equation}
\end{proposition}
\begin{proof}
Follow notation as in the proof of Proposition \ref{normal}. For any fixed $q=(q_1,\ldots,q_n)\in\C[\z]_s^n$, define $w_q=\sum_{i=1}^n\bar z_iq_i$. Moreover, let $\mathbf{q}$ be the coefficient vector of $q$ with respect to the monomial basis $\left[z_1[\z]_s,\ldots,z_n[\z]_s\right]$. 
Let us write $\mathbf{H}_{s}(\y)$ in block form
\begin{equation*}
    \mathbf{H}_{s}(\y)=\begin{bmatrix}
        A&B\\B^*&C
    \end{bmatrix}.
\end{equation*}
Assume that $A\succ0$. Otherwise, use the Moore–Penrose inverse $A^{\dagger}$ instead of $A^{-1}$ in the following. The Schur complement of $A$ is $C-B^*A^{-1}B$, and we have
\begin{align*}
    \mathbf{q}^*(C-B^*A^{-1}B)\mathbf{q}
    &=\langle w_q,w_q\rangle_\y-\left\langle \mathcal{P}_s(w_q),\mathcal{P}_s(w_q)\right\rangle_\y\\
    &=\sum_{1\le i,j\le n}\left(\langle \bar z_iq_i,\bar z_jq_j\rangle_\y-\langle M_i^*(q_i),M_j^*(q_j)\rangle_\y\right)\\
    &=\sum_{1\le i,j\le n}\left(\langle z_jq_i, z_iq_j\rangle_\y-\langle q_i,M_iM_j^*(q_j)\rangle_\y\right)\\
    &=\sum_{1\le i,j\le n}\left(\langle M_j(q_i), M_i(q_j)\rangle_\y-\langle q_i,M_iM_j^*(q_j)\rangle_\y\right)\\
    &=\sum_{1\le i,j\le n}\left(\langle q_i, M_j^*M_i(q_j)\rangle_\y-\langle q_i,M_iM_j^*(q_j)\rangle_\y\right)\\
    &=\sum_{1\le i,j\le n}\left(\langle q_i,[M_j^*,M_i](q_j)\rangle_\y\right).
\end{align*}
We thus conclude that $\mathbf{H}_{s}(\y)\succeq0$ if and only if $\M^{\C}_{s}(\y)\succeq0$ and $\left[[M_j^*,M_i]\right]_{i,j=1}^n\succeq0$. As $M_i$'s are of finite dimension, we have
\begin{equation}
\left[[M_j^*,M_i]\right]_{i,j=1}^n\succeq0\iff[M_i^*,M_j]=0,\,\forall i,j\in[n].
\end{equation}
The desired conclusion then follows from the fact that $M_i$'s are pairwise commuting.
\end{proof}


Let
$b_s(\z)\coloneqq\left[z_1[\z]_s,\ldots,z_n[\z]_s\right]$ which is a subvector of $[\z,\bar{\z}]_{s+1}$. Then,
\begin{equation}\label{joint-normal1}
\mathbf{H}_{s}(\y)=L_{\y}\left(b_s(\z)b_s(\z)^*\right)
\end{equation}
which is a principal submatrix of the $(s+1)$-th order full complex moment matrix $\M^{2\C}_{s+1}(\y)$. For $i\in[n]$, define
\begin{equation}\label{joint-normal2}
\mathbf{H}_{s}(g_i\y)\coloneqq L_{\y}\left(b_s(\z)b_s(\z)^*g_i\right).
\end{equation}
We could strengthen the complex moment relaxation \eqref{cmom} by including the joint normality conditions \eqref{joint-normal1} and \eqref{joint-normal2}:
\begin{equation}\label{ss-cmom}
\tau''_{r,s}\coloneqq\begin{cases}
\inf\limits_{\y}& L_{\y}(f)\\
\rm{s.t.}&\M^{\C}_{r}(\y)\succeq0,\quad y_{\mathbf{0},\mathbf{0}}=1,\\
&\M^{\C}_{r-d^g_i}(g_i\y)\succeq0,\quad i\in[m],\\
&\M^{\C}_{r-d^h_j}(h_j\y)=0,\quad j\in[l],\\
&\mathbf{H}_s(\y)\succeq0,\\
&\mathbf{H}_{s-d^g_i}(g_i\y)\succeq0,\quad i\in[m].
\end{cases}
\end{equation}
By doing so, we again obtain a bilevel hierarchy of strengthened moment relaxations for \eqref{cpop}. In particular, the index $s$ is called the \emph{jointly normal order}. Note that the normality conditions considered in Section~\ref{sec3} correspond to the PSDness of certain principle submatrices of $\mathbf{H}_s(\y)$. Thus, one has $\tau'_{r,s}\le\tau''_{r,s}\le\rho_{\max\,\{r,s+1\}}$.

\begin{remark}
Like the strengthened relaxation considered in Section~\ref{sec3}, one can integrate \eqref{ss-cmom} with correlative sparsity as well as term sparsity to suit for sparse problems.
\end{remark}

Let
$\Sigma_r\coloneqq\left\{\sum_{k}|p_k|^2\,\middle|\, p_k\in\C[\z]_r\right\}$ and
\begin{equation}
\C[\z,\bar\z]_r\coloneqq\left\{p=\sum_{(\b,\g)}p_{\b,\g}\z^{\b}\bar{\z}^{\g}\in\C[\z,\bar\z]\,\middle|\,|\b|,|\g|\le r\right\},
\end{equation}
for $r\in\N$.
Moreover, define
\begin{equation}
\mathcal E\coloneqq\mathbb C[\z]+\sum_{i=1}^n\bar z_i\,\mathbb C[\z]\text{ and }\mathcal E_s\coloneqq\mathbb C[\z]_s+\sum_{i=1}^n\bar z_i\,\mathbb C[\z]_s,
\end{equation}
for $s\in\N$. We also define $\mathcal E^2\coloneqq\left\{\sum_{k}|p_k|^2\,\middle|\, p_k\in\mathcal E\right\}$ and $\mathcal E_s^2\coloneqq\left\{\sum_{k}|p_k|^2\,\middle|\, p_k\in\mathcal E_s\right\}$.
Then the dual of \eqref{ss-cmom} is
\begin{equation}\label{ss-sos}
\begin{cases}
\sup\limits_{\lambda\in\R}& \lambda\\
\,\rm{s.t.}&f-\lambda=\sigma_0+\sum_{i=1}^m\sigma_ig_i+\sum_{j=1}^l\phi_jh_j+\nu_0+\sum_{i=1}^m\nu_ig_i,\\
&\text{where }\sigma_0\in\Sigma_r,\nu_0\in\mathcal E_s^2,\sigma_i\in\Sigma_{r-d^g_i},\\
&\quad\quad\quad\nu_i\in\mathcal E_{s-d^g_i}^2,i\in[m],\,\phi_j\in\C[\z,\bar\z]_{r-d_j^h},j\in[l].
\end{cases}
\end{equation} 

For later use, we first prove the matrix version of Theorem 3.3 in \cite{d2011hermitian}, which might also be of independent interest. Let
\begin{equation}
S\coloneqq\left\{\z\in\mathbb C^n \,\middle|\, h(\z,\bar \z)\coloneqq 1-\|\z\|^2=0\right\}
\end{equation}
be the unit sphere.
\begin{lemma}\label{thm:sphere-matrix-hsos}
Let $H(\z,\bar \z)\in\mathbb C[\z,\bar \z]^{m\times m}$
be an Hermitian polynomial matrix. Assume that
\begin{equation}
    H(\z,\bar \z)\succ 0
    \quad \text{on } S.
\end{equation}
Then there exists a holomorphic polynomial matrix
$P(\z)$
and an Hermitian polynomial matrix
$\Phi(\z,\bar \z)\in \C[\z,\bar \z]^{m\times m}$
such that
\begin{equation}
    H(\z,\bar \z)
    =
    P(\z)^*P(\z)
    +
    \left(1-\|\z\|^2\right)\Phi(\z,\bar \z).
\end{equation}
\end{lemma}

\begin{proof}
Write $H(\z,\bar\z)=\sum_{\a,\b}
    H_{\a,\b}\,\bar \z^\a \z^\b$ with
$H_{\a,\b}\in\C^{m\times m}$ and
$H_{\a,\b}^*=H_{\b,\a}$.
Choose an even integer \(D\) such that
$D\ge \max\,\{|\a|,|\b|\}$
for every pair \((\a,\b)\) appearing in the above expansion.
Introduce one extra complex variable \(t\), and define the bihomogenization
\begin{equation}
    H^{\mathrm{hom}}(\z,t)
    \coloneqq
    \sum_{\a,\b}
    H_{\a,\b}\,
    \bar \z^\a \z^\b
    \bar t^{D-|\a|}
    t^{D-|\b|}.
\end{equation}
Then \(H^{\mathrm{hom}}\) is an \(m\times m\) Hermitian polynomial matrix in
\((\z,t)\), bihomogeneous of bidegree \((D,D)\). In particular, $H^{\mathrm{hom}}(\lambda \z,\lambda t)
=|\lambda|^{2D}H^{\mathrm{hom}}(\z,t),\forall \lambda\in\mathbb C$, and $H^{\mathrm{hom}}(\z,1)=H(\z,\bar \z)$.
Let $R(\z,t)\coloneqq\|\z\|^2-|t|^2$.
We first observe that \(H^{\mathrm{hom}}\) is positive definite on every nonzero point of the cone defined by $R(\z,t)=0$.
Indeed, suppose \(R(\z,t)=0\) and \((\z,t)\neq \mathbf{0}\). Then \(t\neq0\), since otherwise
\(\|\z\|=0\), forcing \((\z,t)=\mathbf{0}\). Hence \(\bw=\z/t\) is well-defined and
we have $\|\bw\|^2
=\frac{\|\z\|^2}{|t|^2}=1$.
For \(t\neq0\), the definition of \(H^{\mathrm{hom}}\) gives
\begin{equation}
    H^{\mathrm{hom}}(\z,t)
    =
    |t|^{2D}
    H\left(\frac \z t,\frac{\bar \z}{\bar t}\right)
    =
    |t|^{2D}H(\bw,\bar \bw).
\end{equation}
Since \(H(\bw,\bar \bw)\succ0\) on the unit sphere and \(|t|^{2D}>0\), it follows that $H^{\mathrm{hom}}(\z,t)\succ0$ on every nonzero point of the cone defined by \(R=0\).
The matrix \(H^{\mathrm{hom}}\) need not be positive definite away from this cone. We now add a large positive term that vanishes on the cone. Since \(D\) is even, $R(\z,t)^D\ge0$
for all \((\z,t)\). For \(C>0\), define
\begin{equation}
F_C(\z,t)
    \coloneqq
    H^{\mathrm{hom}}(\z,t)
    +
    CR(\z,t)^D I_m,
\end{equation}
where $I_m$ is the $m\times m$ identity matrix.
Both terms are bihomogeneous of bidegree \((D,D)\), so \(F_C\) is also Hermitian and bihomogeneous of bidegree \((D,D)\).
We claim that for \(C\) sufficiently large,
$F_C(\z,t)\succ0$ for all $(\z,t)\neq\mathbf{0}$.
By bihomogeneity, it suffices to prove this on the compact sphere
\begin{equation}
    \Omega
    \coloneqq
    \{(\z,t)\in\mathbb C^{n+1}:\|\z\|^2+|t|^2=1\}.
\end{equation}
On the compact set $\Omega_0\coloneqq\{(\z,t)\in\Omega:R(\z,t)=0\}$,
we have already shown that
$H^{\mathrm{hom}}(\z,t)\succ0$.
Therefore, by continuity of the smallest eigenvalue, there exist constants
\(\eta>0\) and \(\delta>0\) such that
$H^{\mathrm{hom}}(\z,t)\succeq \delta I_m$
whenever $(\z,t)\in\Omega$ and $|R(\z,t)|\le\eta$.
On the complementary compact set
$\Omega_1
    \coloneqq
    \{(\z,t)\in\Omega: |R(\z,t)|\ge\eta\},$
the smallest eigenvalue of \(H^{\mathrm{hom}}\) is bounded from below. Hence, there exists
\(\kappa\ge0\) such that
$H^{\mathrm{hom}}(\z,t)\succeq -\kappa I_m$ for all $(\z,t)\in\Omega_1$.
On \(\Omega_1\), since \(D\) is even, $R(\z,t)^D\ge \eta^D$.
Thus, if $C>\frac{\kappa}{\eta^D}$,
then
\begin{equation}
    F_C(\z,t)
    =
    H^{\mathrm{hom}}(\z,t)
    +
    CR(\z,t)^D I_m
    \succ0
\end{equation}
on \(\Omega_1\). On the region \(|R|\le\eta\), the same conclusion already follows from
\(H^{\mathrm{hom}}\succeq\delta I_m\). Hence, for this choice of \(C\),
$F_C(\z,t)\succ0$ for every $(\z,t)\in\Omega$.
By bihomogeneity, the same holds for every nonzero \((\z,t)\in\mathbb C^{n+1}\).

We may now apply the matrix-valued Catlin-D'Angelo's stabilization theorem to the
positive definite bihomogeneous polynomial matrix \(F_C\) \cite{catlin1999isometric,d1999holomorphic}: There exists an
integer \(q\ge0\) and a holomorphic polynomial matrix \(A(\z,t)\) such that
\begin{equation}
    \left(\|\z\|^2+|t|^2\right)^q
    F_C(\z,t)
    =
    A(\z,t)^*A(\z,t).
\end{equation}
Letting \(t=1\) and using \(H^{\mathrm{hom}}(\z,1)=H(\z,\bar \z)\), we obtain
\begin{equation}
    (1+\|\z\|^2)^q
    \left(
        H(\z,\bar \z)
        +
        C(\|\z\|^2-1)^D I_m
    \right)
    =
    A(\z,1)^*A(\z,1).
\end{equation}
Define $P(\z)\coloneqq2^{-q/2}A(\z,1)$. Then
\begin{equation}
    P(\z)^*P(\z)
    =
    2^{-q}(1+\|\z\|^2)^q H(\z,\bar \z)
    +
    2^{-q}C(1+\|\z\|^2)^q(\|\z\|^2-1)^D I_m.
\end{equation}
Subtracting \(H\), we get
\begin{equation}
    P^*P-H=
    \left(2^{-q}(1+\|\z\|^2)^q-1\right)H
    +
    2^{-q}C(1+\|\z\|^2)^q(\|\z\|^2-1)^D I_m.
\end{equation}
The polynomial
$2^{-q}(1+\|\z\|^2)^q-1$
vanishes at \(\|\z\|^2=1\), and hence is divisible by \(\|\z\|^2-1\). Therefore,
$P^*P-H$ is divisible by \(\|\z\|^2-1\) and so 
there exists an Hermitian polynomial matrix \(B(\z,\bar \z)\) such that $P^*P-H=-\left(1-\|z\|^2\right)B(\z,\bar \z)$.
Equivalently, $H=P^*P+\left(1-\|\z\|^2\right)B$.
\end{proof}

We next establish an $\mathcal{E}^2$-certificate for a self-conjugate polynomial to be nonnegative on the unit sphere $S$.

\begin{theorem}\label{thm:sphere-finite-global minimizers}
Let $f\in\C[\z,\bar\z]$ be a self-conjugate polynomial that is nonnegative on the unit sphere $S$. Assume that the zero set of $f$ is nonempty and finite, and that every zero as a global minimizer is nondegenerate on \(S\).
Then there exists a polynomial $\nu\in\mathcal E^2$
and a self-conjugate polynomial \(\psi\) such that
\begin{equation}\label{sec5:eq1}
f=\nu+h\psi.
\end{equation}
\end{theorem}

\begin{proof}
Denote the zero set of $f$ by $Z\coloneqq\{\bu^{(1)},\cdots,\bu^{(a)}\}$, and let
\begin{equation}
I_{\C[\z]}(Z)
\coloneqq
\left\{p\in\mathbb C[\z]\,\middle|\, p(\bu^{(i)})=0,\ i=1,\ldots,a\right\}
\end{equation}
be the holomorphic vanishing ideal of \(Z\). Since \(Z\) is finite, we may choose generators
$g_1,\ldots,g_b\in\mathbb C[\z]$
such that $I_{\C[\z]}(Z)=(g_1,\ldots,g_b)$.
For each point $\bu^{(i)}$, define the holomorphic maximal ideal
\begin{equation}
\mathfrak m_i\coloneqq
(z_1-u_1^{(i)},\ldots,z_n-u_n^{(i)})\subseteq\C[\z].
\end{equation}
Because the points $\{\bu^{(1)},\cdots,\bu^{(a)}\}$ are distinct, the maximal ideals $\{\mathfrak m_1,\cdots,\mathfrak m_a\}$ as well as $\{\mathfrak m_1^2,\cdots,\mathfrak m_a^2\}$ are pairwise comaximal. Therefore, the Chinese remainder theorem gives
\begin{equation}
\C[\z]/\bigcap_{i=1}^a\mathfrak{m}_i^2\cong\prod_{i=1}^b\C[\z]/\mathfrak{m}_i^2.
\end{equation}
It follows that there exist Hermite interpolation polynomials $L_1,\ldots,L_a\in\mathbb C[\z]$
satisfying
\begin{equation}
L_i(\bu^{(j)})=\delta_{ij},
\qquad
dL_i(\bu^{(j)})=0
\quad
\text{for all }i,j.
\end{equation}
Equivalently, we prescribe
\begin{equation}
L_i\equiv1\pmod{\mathfrak m_i^2},
\qquad
L_i\equiv0\pmod{\mathfrak m_j^2}\quad\text{ for }i\ne j.
\end{equation}
Now define
\begin{equation}
s_{i,k}(\z)=L_i(\z)(z_k-u_k^{(i)}),\quad r_{i,k}(\z,\bar \z)=L_i(\z)(\bar z_k-\bar u_k^{(i)}),
\end{equation}
for $i\in[a]$ and $k\in[n]$.
Let $V$ be the column vector formed by
$\{g_j,\ s_{i,k},\ r_{i,k}:
j=1,\ldots,b,\ i=1,\ldots,a,\ k=1,\ldots,n\}$
and write $V=[v_1,\ldots,v_m]^{\intercal}$.
Every component of \(V\) lies in
$\mathcal E.$
We claim that the components of \(V\) generate the defining ideal $I_A(Z)$ of \(Z\) in \(A\coloneqq\C[\z,\bar\z]/(h)\), where
\begin{equation}
I_A(Z)\coloneqq
\left\{p\in A\,\middle|\, p(\bu^{(i)},\bar \bu^{(i)})=0,\ i=1,\ldots,a\right\}.
\end{equation}
Clearly, $(v_1,\ldots,v_m)_{A}\subseteq I_A(Z)$,
because each component of \(V\) vanishes on \(Z\).
Conversely, take any \(p\in I_A(Z)\) and choose a representative, still denoted by \(p\), in \(\C[\z,\bar\z]\). For each \(i\), since \(p(\bu^{(i)},\bar \bu^{(i)})=0\), the polynomial identity
\begin{equation}
p(\z,\bar \z)
=
\sum_{k=1}^n \left(z_k-u_k^{(i)}\right)\xi_{i,k}(\z,\bar \z)
+
\sum_{k=1}^n \left(\bar z_k-\bar u_k^{(i)}\right)\eta_{i,k}(\z,\bar \z)
\end{equation}
holds for suitable \(\xi_{i,k},\eta_{i,k}\in \C[\z,\bar\z]\). This is just Taylor's formula in the polynomial ring \(\C[\z,\bar\z]\) around the point $(\bu^{(i)},\bar \bu^{(i)})$.
Multiplying by \(L_i\) gives
\begin{equation}
L_ip=\sum_k \left(s_{i,k}\xi_{i,k}+r_{i,k}\eta_{i,k}\right)
\in (v_1,\ldots,v_m)_{A}.
\end{equation}
Now put $L=\sum_{i=1}^a L_i$.
Since $L(\bu^{(i)})=1$ for every $i$,
we have $1-L\in I_{\C[\z]}(Z)\subseteq (v_1,\ldots,v_m)_{A}$.
Hence,
\begin{equation}
p=Lp+(1-L)p=
\sum_i L_ip+(1-L)p
\in (v_1,\ldots,v_m)_{A}.
\end{equation}
Thus $I_A(Z)=(v_1,\ldots,v_m)_{A}$.

Let \(\bu\in Z\) be any zero of $f$. Since \(\bu\) is a constrained global minimizer of \(f|_S\), it holds
$d(f|_S)(\bu,\bar\bu)=0.$
Equivalently, the class of \(f\) in the localization ring \(\mathcal O_{S,\bu}\) of $A$ at $\bu$ has zero value and zero first-order part. Thus $f\in \mathfrak m_{S,\bu}^2$ locally at \(\bu\), where
$\mathfrak m_{S,\bu}
\coloneqq
\{p\in\mathcal O_{S,\bu}\mid p(\bu,\bar\bu)=0\}.$
Since \(Z\) is finite,
\begin{equation}
I_A(Z)=\bigcap_{i=1}^a \mathfrak m_{S,\bu^{(i)}}.
\end{equation}
The maximal ideals \(\mathfrak m_{S,\bu^{(i)}}\) are pairwise comaximal, and hence
\begin{equation}
I_A(Z)^2
=
\left(\bigcap_i\mathfrak m_{S,\bu^{(i)}}\right)^2
=
\bigcap_i \mathfrak m_{S,\bu^{(i)}}^2.
\end{equation}
Therefore, $f\in I_A(Z)^2$. It follows that \(f\) is a finite sum of products of two elements in $\{v_1,\ldots,v_m\}$.
Because \(I_A(Z)\) is closed under conjugation, each \(v_i\) is an \(A\)-linear combination of the \(\bar v_j\)'s. Consequently, every product of two generators can be rewritten as a sum of terms of the form $c_{ij}\bar v_iv_j$.
Hence, there exists a polynomial matrix \(G\in A^{m\times m}\) such that
\begin{equation}\label{sec5:eq2}
f=V^*GV
\quad
\text{in }A.
\end{equation}
Since \(f=\bar f\), we may assume $G=G^*$.

Fix any \(\bu\in Z\). Since \(V(\bu)=\mathbf{0}\), for a smooth curve
\begin{equation}
\gamma(t)\subset S,
\qquad
\gamma(0)=\bu,
\qquad
\dot\gamma(0)=\bv\in T_\bu S,
\end{equation}
we have
\begin{equation}
V(\gamma(t))
=t\,dV(\bu)[\bv]+O(t^2),
\end{equation}
and
\begin{equation}
G(\gamma(t))=G(\bu)+O(t).
\end{equation}
Along \(\gamma\), \eqref{sec5:eq2} gives
\begin{equation}
f(\gamma(t))
=
V(\gamma(t))^*G(\gamma(t))V(\gamma(t)).
\end{equation}
Therefore,
\begin{equation}
f(\gamma(t))
=
t^2
\bigl(dV(\bu)[\bv]\bigr)^*
G(\bu)
\bigl(dV(\bu)[\bv]\bigr)
+
O(t^3).
\end{equation}
Thus,
\begin{equation}
\frac12\operatorname{Hess}_{S}f(\bu)[\bv,\bv]
=
\bigl(dV(\bu)[\bv]\bigr)^*
G(\bu)
\bigl(dV(\bu)[\bv]\bigr).
\end{equation}
Since \(\bu\) is a nondegenerate minimizer,
$\operatorname{Hess}_{S}f(\bu)\succ0$ on $T_\bu S$.
Hence $\bw^*G(\bu)\bw>0$ for every $\mathbf{0}\ne \bw\in W_\bu\coloneqq\operatorname{im}\bigl(dV(\bu):T_\bu S\to\mathbb C^m\bigr)$.

A syzygy of \(V\) is a row vector \(U\in A^{m}\) such that $UV=0$ in $A$.
If \(U\) is any matrix whose rows are syzygies, then
$V^*U^*UV=(UV)^*(UV)=0$ in $A$.
Therefore, replacing \(G\) by
$G+U^*U$
does not change the represented polynomial \(V^*GV\).
We need enough syzygies to control the directions transverse to \(W_\bu\). Let
\begin{equation}
W_\bu^\perp
\coloneqq
\{\ell\in(\mathbb C^m)^*\mid\ell(\bw)=0\text{ for all }\bw\in W_\bu\}.
\end{equation}
We claim that every \(\ell\in W_\bu^\perp\) occurs as the value at \(\bu\) of a syzygy.
Indeed, let us write
$\ell=(\ell_1,\ldots,\ell_m)$.
The condition \(\ell\in W_\bu^\perp\) means
$d\left(\sum_{i=1}^m \ell_i v_i\right)(\bu)=0$
on \(T_\bu S\). Since we also have \(v_i(\bu)=0\), the element $q=\sum_i \ell_i v_i$
belongs to \(\mathfrak m_{S,\bu}^2\). But
$\mathfrak m_{S,\bu}=(v_1,\ldots,v_m)_{A}$
locally, and hence $q=\sum_i c_i v_i$
with $c_i\in\mathfrak m_{S,\bu}$.
Therefore,
$\sum_i(\ell_i-c_i)v_i=0$
locally, and the row vector $U_\ell=(\ell_i-c_i)_{i=1}^m$
is a local syzygy satisfying $U_\ell(\bu)=\ell$.
By clearing local denominators and using Hermite interpolation over the finite set \(Z\), these local syzygies can be represented by global polynomial syzygies.

Choose finitely many global syzygy rows and put them into a matrix \(U_Z\) so that, for every \(\bu\in Z\),
$\ker U_Z(\bu)=W_\bu$. Now define
\begin{equation}
G_\mu=G+\mu U_Z^*U_Z.
\end{equation}
For every \(\bu\in Z\), the matrix \(G(\bu)\) is positive definite on \(W_\bu\), while
$U_Z(\bu)^*U_Z(\bu)$
is positive definite on directions transverse to \(W_\bu\). Therefore, for sufficiently large \(\mu>0\), $G_\mu(\bu)\succ0$ for all $\bu\in Z$.
Since
$V^*G_\mu V=V^*GV$ in $A$,
replacing \(G\) by \(G_\mu\), we may henceforth assume
$G(\bu)\succ0$ for all $\bu\in Z$.
By continuity, there exists a neighborhood \(O\) of \(Z\) in \(S\) such that
$G(\z,\bar \z)\succ0$ for all $\z\in O$.

On \(S\setminus O\), we have $f=V^*GV>0$ because \(f\) vanishes exactly on \(Z\). 
Consider the Hermitian matrix
\begin{equation}
C(\z,\bar \z)
=
\|V(\z,\bar \z)\|^2I_m
-
V(\z,\bar \z)V(\z,\bar \z)^*,
\end{equation}
where $I_m$ is the $m\times m$ identity matrix.
For every vector \(\bv\in\mathbb C^m\),
$\bv^*C\bv=\|V\|^2\|\bv\|^2-|V^*\bv|^2\ge0$
by Cauchy's inequality. Hence $C\succeq0$.
Moreover, $CV=\mathbf{0}$,
and therefore, $V^*CV=0$.
Define $G_\lambda=G+\lambda C$.
It holds $V^*G_\lambda V=V^*A_0V$ in $A$.
We claim that for sufficiently large \(\lambda>0\),
\begin{equation}
G_\lambda(\z,\bar \z)\succ0,
\qquad
\forall \z\in S\setminus O.
\end{equation}
Fix \(\z\in S\setminus O\) and write $\bv=V(\z,\bar \z)\ne\mathbf{0}$.
Since $\bv^*G(\z,\bar \z)\bv=f(\z,\bar \z)>0$,
\(G(\z,\bar \z)\) is positive definite on the line \(\mathbb C\bv\). The correction $\|\bv\|^2I_m-\bv\bv^*$
vanishes on \(\mathbb C\bv\) and is positive definite on \(\bv^\perp\). Therefore, for this fixed \(\z\), \(G_\lambda(\z,\bar \z)=G(\z,\bar \z)+\lambda(\|\bv\|^2I_m-\bv\bv^*)\) is positive definite for all sufficiently large \(\lambda\). Furthermore,
because \(S\setminus O\) is compact
and the entries of \(G\) are continuous, one single sufficiently large \(\lambda\) works for all \(\z\in S\setminus O\).
On \(O\), \(G\succ0\), and \(C\succeq0\), so we have $G_\lambda\succ0$ on $O$.
Thus, for sufficiently large \(\lambda\),
\begin{equation}
G_\lambda(\z,\bar \z)\succ0,
\qquad
\forall \z\in S.
\end{equation}
Replacing \(G\) by \(G_\lambda\), we obtain a Hermitian polynomial matrix \(G\) such that
\begin{equation}
G(\z,\bar \z)\succ0,\qquad\forall \z\in S,
\end{equation}
and
\begin{equation}
f=V^*GV
\quad
\text{in }A.
\end{equation}
Equivalently, 
$f=V^*GV+h\phi$ for some self-conjugate polynomial \(\phi\in \C[\z,\bar\z]\).

Now applying Lemma~\ref{thm:sphere-matrix-hsos} to the matrix \(G\), we get
$G=P^*P+h\Phi$ for some holomorphic polynomial matrix $P(\z)$ and some Hermitian polynomial matrix
$\Phi(\z,\bar \z)$. Therefore,
\begin{equation}
f=
V^*P^*PV+hV^*\Phi V+h\phi=
\|PV\|^2+h\psi,
\end{equation}
where $\psi=V^*\Phi V+\phi$ is self-conjugate.
Finally, each component of \(PV\) lies in \(\mathcal E\) because: $P(\z)$ has holomorphic entries, each component of \(V\) lies in \(\mathcal E\), and \(\mathcal E\) is closed under addition and multiplication by holomorphic polynomials.
Therefore, writing the components of \(PV\) as
$q_k(\z,\bar \z)
=
p_{k,0}(\z)+\sum_{i=1}^n\bar z_i\,p_{k,i}(\z)$,
we obtain
$f=\sum_k |q_k|^2+h\psi$ as desired.
\end{proof}

Now we are in a position to prove the finite convergence result for optimizing a self-conjugate polynomial over spheres.
\begin{theorem}
Let $f\in\C[\z,\bar\z]$ be a self-conjugate polynomial and consider the minimization of $f$ over the unit sphere $S$: $\min_{\z\in S}f(\z,\bar\z)$.
Assume that the set of global minimizers
is nonempty and finite, and that every global minimizer is nondegenerate on \(S\). Then the strengthened complex moment hierarchy \eqref{ss-cmom} has finite convergence, that is, $\tau''_{r,s}=f_{\min}$ holds true for any $r\in\N$ and sufficiently large $s$. Furthermore, for any $r\in\N$ and sufficiently large $s$, there exists an integer \(t\le s\) such that, every optimal solution \(\y\) of \eqref{ss-cmom} satisfies
\begin{equation}
\rank\,\M_{t}^{\mathbb C}(\y)=\rank\,\M_{t-2}^{\mathbb C}(\y).
\end{equation}
In other words, the hierarchy admits flat truncation.
\end{theorem}

\begin{proof}
Set $\tilde{f}=f-f_{\min}$.
By Theorem~\ref{thm:sphere-finite-global minimizers}, there exists $\nu\in\mathcal E^2$ and a self-conjugate polynomial $\psi$ such that $\tilde{f}=\nu+h\psi$. 
Choose \(s\) large enough so that all terms in \(\nu+h\psi\) are representable in \eqref{ss-sos}. Hence, by weak duality, $\tau''_{r,s}=f_{\min}$. 

Denote the set of global minimizers by $Z\coloneqq\{\bu^{(1)},\cdots,\bu^{(a)}\}$, and let
$I_{\C[\z]}(Z)$
be the holomorphic vanishing ideal of \(Z\).
Fix any holomorphic polynomial
$p\in I_{\C[\z]}(Z)$.
We claim that, for sufficiently large $s$, every optimal solution \(\y\) of \eqref{ss-cmom} satisfies
$p\in\ker \M_{\deg(p)}^{\mathbb C}(\y)$.
To prove this, we first construct a nonnegative domination certificate. Because \(p\) vanishes on \(Z\), near each global minimizer \(\bu^{(i)}\) we have
\begin{equation}
|p(\z)|^2=O(\operatorname{dist}_S(\z,\bu^{(i)})^2).
\end{equation}
Since \(\bu^{(i)}\) is a nondegenerate global minimizer of \(f|_S\),
\begin{equation}
\tilde{f}(\z,\bar\z)\ge c_i\operatorname{dist}_S(\z,\bu^{(i)})^2
\end{equation}
near \(\bu^{(i)}\), for some \(c_i>0\).
Away from a small neighborhood of \(Z\), compactness gives $\tilde{f}(\z,\bar\z)\ge \delta>0$.
Therefore, choosing \(C_p>0\) to be sufficiently large, we obtain
$\omega_p\coloneqq C_p\tilde{f}-|p|^2\ge0$ on $S$,
and the zero set of \(\omega_p\) on \(S\) is exactly \(Z\).
Moreover, \(C_p\) can be chosen large enough so that at every \(\bu^{(i)}\),
$\operatorname{Hess}_S \omega_p(\bu^{(i)})\succ0$.
Indeed,
\begin{equation}
\operatorname{Hess}_S \omega_p(\bu^{(i)})
=
C_p\operatorname{Hess}_S \tilde{f}(\bu^{(i)})
-
\operatorname{Hess}_S |p|^2(\bu^{(i)}),
\end{equation}
and the first term is positive definite for sufficiently large \(C_p\).
Thus, \(\omega_p\) satisfies the same properties as \(\tilde{f}\): it is nonnegative on \(S\), has the finite zero set \(Z\), and has a positive definite Hessian sphere at every zero. Applying Theorem~\ref{thm:sphere-finite-global minimizers}, we get
\begin{equation}\label{sec5:eq3}
C_p\tilde{f}-|p|^2
=\nu^{(p)}
+
h\psi_p,
\end{equation}
where $\nu^{(p)}
\in\mathcal E^2.$
Choose \(s\) large enough so that \(\nu^{(p)}+h\psi_p\) is representable in \eqref{ss-sos}.
Now let \(\y\) be an optimal solution of \eqref{ss-cmom} at such a level. Applying \(L_\y\) to \eqref{sec5:eq3} and using $L_\y(h\psi_p)=0$, we have
\begin{equation}\label{sec:eq4}
C_pL_\y(\tilde{f})-L_\y(|p|^2)
=L_\y(\nu^{(p)}).
\end{equation}
Since $L_\y(\tilde{f})=0$ and
$L_\y(\nu^{(p)})\ge0$, \eqref{sec:eq4} gives
$-L_\y(|p|^2)\ge0$.
But \(\M_{\deg(p)}^{\mathbb C}(\y)\succeq0\), so
$L_\y(|p|^2)\ge0$.
It follows $L_\y(|p|^2)=0$, which implies that the coefficient vector of \(p\) lies in the kernel of the moment matrix $\M_{\deg(p)}^{\mathbb C}(\y)$.
This proves the claim.

The ideal \(I_{\C[\z]}(Z)\subset\mathbb C[\z]\) is zero-dimensional and radical. Let $\Gamma$ be the reduced Gr\"obner basis
of \(I_{\C[\z]}(Z)\) with respect to some total-degree order, and let $\mathcal B$ be the corresponding standard monomial basis of the quotient ring $\mathbb C[\z]/I_{\C[\z]}(Z)$.
Fix $t=\max_{b\in\mathcal B}\deg(b)+2$.
For every monomial \(\z^\a\) with \(|\a|\le t\), polynomial division by the Gröbner basis $\Gamma$ gives
$\z^\a-\eta_\a(\z)\in I_{\C[\z]}(Z)$,
where
$\eta_\a\in\operatorname{span}\mathcal B$ is the normal form of \(\z^\a\).
Thus, $\z^\a-\eta_\a(\z)\in\ker \M_{t}^{\mathbb C}(\y)$ for every \(|\a|\le t\).
It follows that the column of \(\M_{t}^{\mathbb C}(\y)\) indexed by \(\z^\a\) is a linear combination of the columns indexed by the standard monomials in \(\mathcal B\).
Hence, $\operatorname{rank}\M_{t}^{\mathbb C}(\y)=\operatorname{rank}\M_{t-2}^{\mathbb C}(\y)$.
\end{proof}

\begin{remark}
The nonemptyness, finiteness, and nondegeneracy of global minimizers are generic properties in the space of coefficients of $f$. Thus, both finite convergence and flat truncation generically occur for the strengthened complex moment hierarchy \eqref{ss-cmom} of complex polynomial optimization over spheres.
\end{remark}

Before closing this section, we compare the maximal PSD block sizes of different relaxations in Table~\ref{size}, where `Mom' stands for the complex moment relaxation \eqref{cmom}, `S-Mom' stands for the first strengthened complex moment relaxation \eqref{s-cmom}, `SS-Mom' stands for the second strengthened complex moment relaxation \eqref{ss-cmom}, and `F-Mom' stands for the full complex moment relaxation \eqref{fcmom}. We could see that, when $s$ is chosen to be small, the block size is dominated by the moment matrix constraint. As a result, the maximal block sizes of two strengthened relaxations are close to that of the pruned complex moment relaxation and are much smaller than that of the full complex moment relaxation.

\begin{table}[htbp]\label{size}
\caption{The maximal PSD block sizes of different relaxations.}
    \centering
    \renewcommand\arraystretch{1.8}
    \begin{tabular}{|c|c|c|c|}
    \hline
         Mom& S-Mom& SS-Mom & F-Mom\\
         \hline
       $\binom{n+r}{r}$&$\max\left\{\binom{n+r}{r}, 2\binom{n+s}{s}\right\}$&$\max\left\{\binom{n+r}{r}, (n+1)\binom{n+s}{s}\right\}$&$\binom{2n+r}{r}$\\
       \hline
    \end{tabular}
\end{table}

\section{Numerical experiments}
The strengthened complex moment hierarchies have been implemented in the Julia package {\tt TSSOS}\footnote{{\tt TSSOS} is freely available at \href{https://github.com/wangjie212/TSSOS}{https://github.com/wangjie212/TSSOS}.}. In this section, we evaluate their performance on diverse complex polynomial optimization problems using {\tt TSSOS} where {\tt Mosek} 10.2 \cite{mosek} is employed as an SDP backend with default settings. When presenting the results, the column labelled by `opt' records optima of SDP relaxations and the column labelled by `time' records computational time in seconds. Lower bounds are marked by `*' if global optimality is certified through the flatness condition.
Moreover, the symbol `-' indicates that {\tt Mosek} runs out of memory. Unless otherwise stated, the (jointly) normal order $s$ is set to 1.
All numerical experiments were performed on a desktop computer with an Intel(R) Core(TM) i9-10900 CPU@2.80GHz and 64G RAM.

\subsection{Minimizing a random complex quartic polynomial with unit-norm variables} 
Let us minimize a random complex quartic polynomial with unit-norm variables:
\begin{equation}
\begin{cases}
\inf\limits_{\z\in\C^{n}} &[\z]_2^{*}Q[\z]_2\\
\,\,\,\rm{s.t.}&|z_i|^2=1,\quad i=1,\ldots,n,
\end{cases}
\end{equation}
where $Q\in\C^{|[\z]_2|\times|[\z]_2|}$ ($|[\z]_2|$ is the cardinality of $[\z]_2$) is a random Hermitian matrix whose entries (both real and imaginary parts) are selected with respect to the uniform probability distribution on $[0,1]$. For each $n\in\{5,10,15\}$, we solve three random instances and present the results in Table \ref{tab3}. By Proposition \ref{unitnorm}, the second-order strengthening relaxation \eqref{s-cmom} with $s=1$ yields the same lower bound with the second-order pruned relaxation and we thus ignore the results of `S-Mom' for this example.

The table shows that the joint-normality strengthening substantially improves the second-order pruned relaxation. For $n=5$ and $n=10$, `SS-Mom' attains the same certified optimum as the third-order pruned and full complex relaxations, while requiring markedly less computational time. For $n=15$, the third-order relaxation runs out of memory, whereas `SS-Mom' remains tractable; it matches the full relaxation in two instances and gives a substantially stronger bound than the second-order pruned relaxation in the remaining instance. Thus, `SS-Mom' recovers most or all of the strength of the more expensive relaxations at significantly lower cost.


\begin{table}[htbp]\label{tab3}
\caption{Minimizing a random complex quartic polynomial with unit-norm variables.}
\renewcommand\arraystretch{1.2}
\centering
\resizebox{\linewidth}{!}{
\begin{tabular}{c|c|c|c|c|c|c|c|c|c}
\multirow{2}*{$n$}&\multirow{2}*{trial}&\multicolumn{2}{c|}{Mom ($r=2$)}&\multicolumn{2}{c|}{Mom ($r=3$)}&\multicolumn{2}{c|}{SS-Mom ($r=2$)}&\multicolumn{2}{c}{F-Mom ($r=2$)}\\
\cline{3-10}
&&opt&time&opt&time&opt&time&opt&time\\
\hline
\multirow{3}*{$5$}&1&-24.545&0.03&-22.290*&0.20&-22.290*&0.05&-22.290*&0.13\\
&2&-21.750&0.02&-19.102*&0.18&-19.102*&0.05&-19.102*&0.14\\
&3&-20.838&0.02&-19.281*&0.25&-19.281*&0.04&-19.281*&0.19\\
\hline
\multirow{3}*{$10$}&1&-151.45&1.14&-111.78*&1020&-111.78*&2.11&-111.78*&13.8\\
&2&-149.67&1.21&-93.281*&1121&-93.281*&2.23&-93.281*&13.9\\
&3&-152.79&1.19&-98.605*&1002&-98.605*&1.87&-98.605*&10.2\\
\hline
\multirow{3}*{$15$}&1&-511.33&42.3&-&-&-311.51*&54.7&-311.51*&1085\\
&2&-513.23&46.4&-&-&-290.47&59.3&-280.87*&1019\\
&3&-515.47&47.4&-&-&-302.49*&72.6&-302.49*&1004\\
 \end{tabular}}
\end{table}


\subsection{Minimizing a random complex quartic polynomial on a sphere}
Let us minimize a random complex quartic polynomial on a unit sphere:
\begin{equation}
\begin{cases}
\inf\limits_{\z\in\C^{n}} &[\z]_2^{*}Q[\z]_2\\
\,\,\,\rm{s.t.}&|z_1|^2+\cdots+|z_n|^2=1,
\end{cases}
\end{equation}
where $Q\in\C^{|[\z]_2|\times|[\z]_2|}$ is a random Hermitian matrix whose entries (both real and imaginary parts) are selected with respect to the uniform probability distribution on $[0,1]$. For each $n\in\{5,10,15\}$, we solve three random instances and present the results in Table \ref{tab4}. 

The table shows that both strengthenings improve the second-order pruned relaxation, with `SS-Mom' consistently producing the strongest bounds. For $n=5$ and $n=10$, `SS-Mom' matches the full complex relaxation and certifies global optimality for all instances, while being substantially faster for $n=10$. For $n=15$, both the third-order pruned and full complex relaxations run out of memory, whereas `SS-Mom' remains tractable and certifies global optimality in all three instances. These results demonstrate that the joint-normality strengthening can recover the strength of the full relaxation with markedly better scalability.


\begin{table}[htbp]\label{tab4}
\caption{Minimizing a random complex quartic polynomial on a unit sphere.}
\renewcommand\arraystretch{1.2}
\centering
\resizebox{\linewidth}{!}{
\begin{tabular}{c|c|c|c|c|c|c|c|c|c|c|c}
\multirow{2}*{$n$}&\multirow{2}*{trial}&\multicolumn{2}{c|}{Mom ($r=2$)}&\multicolumn{2}{c|}{Mom ($r=3$)}&\multicolumn{2}{c|}{S-Mom ($r=2$)}&\multicolumn{2}{c|}{SS-Mom ($r=2$)}&\multicolumn{2}{c}{F-Mom ($r=2$)}\\
\cline{3-12}
&&opt&time&opt&time&opt&time&opt&time&opt&time\\
\hline
\multirow{3}*{$5$}&1&-3.8882&0.03&-3.2848&1.16&-3.0001&0.04&-2.8994*&0.06&-2.8994*&0.10\\
&2&-2.7110&0.03&-2.2546&1.25&-2.0609*&0.05&-2.0609*&0.06&-2.0609*&0.10\\
&3&-3.5725&0.04&-3.1123&1.27&-2.8213*&0.04&-2.8213*&0.06&-2.8213*&0.09\\
\hline
\multirow{3}*{$10$}&1&-6.1871&2.10&-&-&-4.8130&2.99&-4.4058*&2.59&-4.4058*&14.3\\
&2&-5.6484&2.23&-&-&-4.1882&2.36&-3.4838*&2.76&-3.4838*&14.2\\
&3&-5.6759&2.24&-&-&-4.3203&2.54&-3.8306*&2.45&-3.8306*&15.9\\
\hline
\multirow{3}*{$15$}&1&-8.4667&57.8&-&-&-6.5468&94.6&-5.0064*&98.5&-&-\\
&2&-7.7119&71.1&-&-&-5.8359&94.1&-4.2665*&104&-&-\\
&3&-7.7191&57.9&-&-&-5.8155&85.4&-4.1569*&112&-&-\\
\end{tabular}}
\end{table}



\subsection{Minimizing a random complex quartic polynomial with correlative sparsity on multi-spheres}
Let us minimize a random complex quartic polynomial with correlative sparsity on multi-spheres:
\begin{equation}
\begin{cases}
\inf\limits_{\z\in\C^{n}} &\sum_{i=1}^l[\z_i]_2^{*}Q_i[\z_i]_2\\
\,\,\,\rm{s.t.}&\|\z_i\|^2=1,\quad i\in[l],
\end{cases}
\end{equation}
where $n=4l+2$, $\z_i\coloneqq\{z_{4i-3},\ldots,z_{4i+2}\}$, and $Q_i\in\C^{|[\z_i]_2|\times|[\z_i]_2|}$ is a random Hermitian matrix whose entries (both real and imaginary parts) are selected with respect to the uniform probability distribution on $[0,1]$. For each $l\in\{5,10,15,20\}$, we solve three random instances and present the results in Table \ref{tab5}. 

The table shows that both sparse strengthenings markedly improve the second-order pruned relaxation at modest additional cost. `SS-Mom' consistently produces much stronger bounds than `Mom' with $r=3$, while being roughly two orders of magnitude faster. It matches the full complex relaxation and certifies global optimality in ten of the twelve instances; in the remaining two instances, its bounds are nearly identical to those of `F-Mom'. Overall, the joint-normality strengthening effectively exploits correlative sparsity and recovers almost all the strength of the full relaxation at substantially lower computational cost.


\begin{table}[htbp]\label{tab5}
\caption{Minimizing a random complex quartic polynomial on multi-spheres.}
\renewcommand\arraystretch{1.2}
\centering
\resizebox{\linewidth}{!}{
\begin{tabular}{c|c|c|c|c|c|c|c|c|c|c|c}
\multirow{2}*{$n$}&\multirow{2}*{trial}&\multicolumn{2}{c|}{Mom ($r=2$)}&\multicolumn{2}{c|}{Mom ($r=3$)}&\multicolumn{2}{c|}{S-Mom ($r=2$)}&\multicolumn{2}{c|}{SS-Mom ($r=2$)}&\multicolumn{2}{c}{F-Mom ($r=2$)}\\
\cline{3-12}
&&opt&time&opt&time&opt&time&opt&time&opt&time\\
\hline
\multirow{3}*{$22$}&1&-18.018&0.34&-14.338&53.7&-14.452&0.44&-12.487*&0.48&-12.487*&2.07\\
&2&-15.922&0.32&-12.635&65.8&-12.851&0.47&-11.390*&0.54&-11.390*&2.09\\
&3&-16.223&0.36&-12.947&54.2&-12.697&0.47&-11.320*&0.55&-11.320*&2.43\\
\hline
\multirow{3}*{$42$}&1&-35.188&0.81&-28.218&108&-28.149&1.03&-24.708*&1.17&-24.708*&4.20\\
&2&-32.871&0.78&-26.513&112&-26.405&1.12&-23.541*&1.42&-23.541*&4.06\\
&3&-33.145&0.70&-27.368&109&-26.410&1.29&-24.456*&1.39&-24.456*&4.29\\
\hline
\multirow{3}*{$62$}&1&-50.784&0.91&-40.280&169&-40.576&1.48&-35.056&1.78&-35.046*&7.22\\
&2&-48.489&0.97&-38.418&149&-38.590&1.63&-33.886*&1.92&-33.886*&5.71\\
&3&-48.370&0.90&-39.454&159&-38.363&1.48&-35.346*&1.61&-35.346*&5.81\\
\hline
\multirow{3}*{$82$}&1&-66.503&1.39&-52.537&204&-52.779&2.05&-46.072&2.69&-46.063*&11.6\\
&2&-64.849&1.41&-51.223&211&-51.518&2.06&-45.242*&2.44&-45.242*&10.6\\
&3&-65.615&1.23&-52.996&239&-51.917&1.83&-46.982*&2.26&-46.982*&9.76\\
\end{tabular}}
\end{table}

\subsection{Smale’s Mean Value conjecture}
The following complex polynomial optimization problem is borrowed from \cite{wang2025real}:
\begin{equation}\label{smale:eq2}
\begin{cases}
\sup\limits_{(\z,u)\in\C^{n+1}} &|u|\\
\,\,\quad\rm{s.t.}&|H(z_i)|\ge|u|,\quad i=1,\ldots,n,\\
&z_1\cdots z_n=\frac{(-1)^n}{n+1},\\
&|z_1|^2+|z_2|^2+\cdots+|z_n|^2=n\left(\frac{1}{n+1}\right)^{\frac{2}{n}},
\end{cases}
\end{equation}
where $H(y)\coloneqq\frac{1}{y}\int_{0}^{y}p(z)\,\mathrm{d}z$ and $p(z)\coloneqq(n+1)(z-z_1)\cdots(z-z_n)$ with $p(0)=1$.
This problem is used in \cite{wang2025real} to verify Smale’s Mean Value conjecture \cite{smale1981fundamental,smale1998mathematical} which is open for $n\ge4$ since 1981. The optimum of \eqref{smale:eq2} is conjectured to be $\frac{n}{n+1}$. We refer the reader to \cite{wang2025real} for more details. Here, we tackle \eqref{smale:eq2} for $n=4$ using the complex moment hierarchy and the strengthened hierarchy \eqref{s-cmom}. The computational results are presented in Table \ref{tab6}, from which we see that the strengthening enables us to achieve global optimality at a much lower relaxation order so that the computational time is dramatically reduced.

\begin{table}[htbp]\label{tab6}
\caption{The results for \eqref{smale:eq2} with $n=4$.}
\renewcommand\arraystretch{1.2}
\centering
\begin{tabular}{|c|c|c|c|c|c|c|}
\hline
\multirow{3}*{Mom}&\multicolumn{2}{c|}{$r=4$}&\multicolumn{2}{c|}{$r=6$}&\multicolumn{2}{c|}{$r=8$}\\
\cline{2-7}
&opt&time&opt&time&opt&time\\
\cline{2-7}
&1.4218&0.12&0.8404&12.3&0.8003&1296\\
\hline
\multirow{3}*{S-Mom}&\multicolumn{2}{c|}{$r=4,s=1$}&\multicolumn{2}{c|}{$r=4,s=2$}&\multicolumn{2}{c|}{$r=4,s=3$}\\
\cline{2-7}
&opt&time&opt&time&opt&time\\
\cline{2-7}
&1.4218&0.13&1.2727&0.15&0.8000&0.27\\
\hline
\end{tabular}
\end{table}

\subsection{The Mordell inequality conjecture}
Our next example concerns the Mordell inequality conjecture due to Birch in 1958: if the numbers $z_1,\ldots,z_n\in\C$ satisfies $|z_1|^2+\cdots+|z_n|^2=n$, then the maximum of $\prod_{1\le i<j\le n}|z_i-z_j|^2$ is $n^n$. This conjecture was proved for $n\le4$ and disproved for $n\ge6$, and so the only remaining open case is when $n=5$. The reader is referred to \cite{wang2025real} for more details.
Without loss of generality, we may eliminate one variable and reformulate the conjecture as the following complex polynomial optimization problem:
\begin{equation}\label{mic}
\begin{cases}
\sup\limits_{\z\in\C^{n-1}} &\prod_{1\le i<j\le n-1}|z_i-z_j|^2\prod_{i=1}^{n-1}|z_i+z_1+\ldots+z_{n-1}|^2\\
\,\,\,\,\,\rm{s.t.}&|z_1|^2+\cdots+|z_{n-1}|^2+|z_1+\ldots+z_{n-1}|^2=n.\\
\end{cases}
\end{equation}
Here, we tackle \eqref{mic} for $n=3,4$ using the complex moment hierarchy and the strengthened hierarchy \eqref{s-cmom}. The computational results are presented in Tables \ref{tab7} and \ref{tab8}, respectively. From the tables, again we see that the strengthening enables us to achieve global optimality
at much lower relaxation orders so that the computational time is dramatically reduced.

\begin{table}[htbp]\label{tab7}
\caption{The results for \eqref{mic} with $n=3$.}
\renewcommand\arraystretch{1.2}
\centering
\begin{tabular}{|c|c|c|c|c|c|c|}
\hline
\multirow{3}*{Mom}&\multicolumn{2}{c|}{$r=10$}&\multicolumn{2}{c|}{$r=14$}&\multicolumn{2}{c|}{$r=18$}\\
\cline{2-7}
&opt&time&opt&time&opt&time\\
\cline{2-7}
&27.347&0.04&27.144&0.17&27.085&0.32\\
\hline
\multirow{3}*{S-Mom}&\multicolumn{2}{c|}{$r=3,s=0$}&\multicolumn{2}{c|}{$r=3,s=1$}&\multicolumn{2}{c|}{$r=3,s=2$}\\
\cline{2-7}
&opt&time&opt&time&opt&time\\
\cline{2-7}
&54.000&0.004&54.000&0.005&27.000&0.005\\
\hline
\end{tabular}
\end{table}

\begin{table}[htbp]\label{tab8}
\caption{The results for \eqref{mic} with $n=4$.}
\renewcommand\arraystretch{1.2}
\centering
\resizebox{\linewidth}{!}{
\begin{tabular}{|c|c|c|c|c|c|c|c|c|c|c|}
\hline
\multirow{3}*{Mom}&\multicolumn{2}{c|}{$r=10$}&\multicolumn{2}{c|}{$r=12$}&\multicolumn{2}{c|}{$r=14$}&\multicolumn{2}{c|}{$r=16$}&\multicolumn{2}{c|}{$r=18$}\\
\cline{2-11}
&opt&time&opt&time&opt&time&opt&time&opt&time\\
\cline{2-11}
&343.67&6.42&326.85&57.8&292.82&205&277.50&669&-&-\\
\hline
\multirow{3}*{S-Mom}&\multicolumn{2}{c|}{$r=6,s=1$}&\multicolumn{2}{c|}{$r=6,s=2$}&\multicolumn{2}{c|}{$r=6,s=3$}&\multicolumn{2}{c|}{$r=6,s=4$}&\multicolumn{2}{c|}{$r=6,s=5$}\\
\cline{2-11}
&opt&time&opt&time&opt&time&opt&time&opt&time\\
\cline{2-11}
&1638.4&0.13&1337.6&0.19&932.20&0.14&582.86&0.15&256.00&0.16\\
\hline
\end{tabular}}
\end{table}

\section{Conclusions}
In this paper, we have shown how to strengthen complex moment relaxations for complex polynomial optimization by employing the normality conditions and the joint normality conditions.
Numerical experiments demonstrate the superior performance of our approach in practice.

On the theoretical side, we have proved that when strengthened by the joint normality condition, the complex moment hierarchy \eqref{ss-cmom} enjoys generic finite convergence and flat truncation for optimizing a self-conjugate polynomial over spheres. At present, the author does not know whether this conclusion could be generalized to complex polynomial optimization problems with general inequality and equality constraints. We thereby formulate it as a conjecture for future studies. 
\begin{conjecture}
Consider \eqref{cpop}. The strengthened complex moment hierarchy \eqref{ss-cmom} has finite convergence and flat truncation in the generic case under the Archimedean condition.
\end{conjecture}
To prove the above conjecture, we need a complex analogy of Marshall's theorem \cite[Theorem 9.5.3]{marshall2008positive}.


\section*{Acknowledgements}
The author thanks Didier Henrion for helpful discussions.
The author is grateful for the comments and suggestions of
the anonymous referees that greatly improved the paper.
The author acknowledges the use of AI for assistance with brainstorming ideas, mathematical development, and drafting the manuscript. The final content, analysis, and conclusions remain the sole responsibility of the
author.

\section*{Funding}
This work was jointly funded by National Key R\&D Program of China under grant No. 2023YFA1009401, Natural Science Foundation of China under grant No. 12201618 and 12171324.

\section*{Conflict of interest}
The authors declare that they have no conflict of interest.

\section*{Data availability}
The author confirms that all data generated or analysed during this study are included in this article.

\bibliographystyle{siamplain}
\bibliography{refer}

\begin{thebibliography}{10}

\bibitem{mosek}
{\sc E.~D. Andersen and K.~D. Andersen}, {\em {The Mosek Interior Point Optimizer for Linear Programming: An Implementation of the Homogeneous Algorithm}}, in High Performance Optimization, vol.~33 of Applied Optimization, Springer US, 2000, pp.~197--232, \url{https://doi.org/10.1007/978-1-4757-3216-0_8}.

\bibitem{aubry2013ambiguity}
{\sc A.~Aubry, A.~De~Maio, B.~Jiang, and S.~Zhang}, {\em Ambiguity function shaping for cognitive radar via complex quartic optimization}, IEEE Transactions on Signal Processing, 61 (2013), pp.~5603--5619.

\bibitem{bienstock2020}
{\sc D.~Bienstock, M.~Escobar, C.~Gentile, and L.~Liberti}, {\em Mathematical programming formulations for the alternating current optimal power flow problem}, 4OR, 18 (2020), pp.~249--292.

\bibitem{catlin1999isometric}
{\sc D.~W. Catlin and J.~P. D'Angelo}, {\em An isometric imbedding theorem for holomorphic bundles}, Mathematical Research Letters, 6 (1999), pp.~43--60.

\bibitem{d2011hermitian}
{\sc J.~P. D'Angelo}, {\em Hermitian analogues of hilbert's 17-th problem}, Advances in Mathematics, 226 (2011), pp.~4607--4637.

\bibitem{dumitrescu2007positive}
{\sc B.~Dumitrescu}, {\em Positive trigonometric polynomials and signal processing applications}, vol.~103, Springer, 2007.

\bibitem{d1999holomorphic}
{\sc J.~P. D’Angelo}, {\em Holomorphic factorization of matrices of polynomials}, in Reproducing Kernels and their Applications, Springer, 1999, pp.~9--19.

\bibitem{fang2021sum}
{\sc K.~Fang and H.~Fawzi}, {\em The sum-of-squares hierarchy on the sphere and applications in quantum information theory}, Mathematical Programming, 190 (2021), pp.~331--360.

\bibitem{fogel2016phase}
{\sc F.~Fogel, I.~Waldspurger, and A.~d’Aspremont}, {\em Phase retrieval for imaging problems}, Mathematical Programming Computation, 8 (2016), pp.~311--335.

\bibitem{goemans2001approximation}
{\sc M.~X. Goemans and D.~Williamson}, {\em Approximation algorithms for max-3-cut and other problems via complex semidefinite programming}, in Proceedings of the thirty-third annual ACM symposium on Theory of computing, 2001, pp.~443--452.

\bibitem{gribling2022bounding}
{\sc S.~Gribling, M.~Laurent, and A.~Steenkamp}, {\em Bounding the separable rank via polynomial optimization}, Linear Algebra and its Applications, 648 (2022), pp.~1--55.

\bibitem{josz2018lasserre}
{\sc C.~Josz and D.~K. Molzahn}, {\em Lasserre hierarchy for large scale polynomial optimization in real and complex variables}, SIAM Journal on Optimization, 28 (2018), pp.~1017--1048.

\bibitem{klep2018minimizer}
{\sc I.~Klep, J.~Povh, and J.~Volcic}, {\em Minimizer extraction in polynomial optimization is robust}, SIAM Journal on Optimization, 28 (2018), pp.~3177--3207.

\bibitem{Las01}
{\sc J.-B. Lasserre}, {\em {Global Optimization with Polynomials and the Problem of Moments}}, SIAM Journal on Optimization, 11 (2001), pp.~796--817.

\bibitem{lasserre2008semidefinite}
{\sc J.~B. Lasserre, M.~Laurent, and P.~Rostalski}, {\em Semidefinite characterization and computation of zero-dimensional real radical ideals}, Foundations of Computational Mathematics, 8 (2008), pp.~607--647.

\bibitem{le2025flat}
{\sc T.~H. Le and M.~T. Ho}, {\em Flat extension technique for moment matrices of positive linear functionals over mixed polynomials and an application in quantum information}, Annals of Functional Analysis, 16 (2025), p.~30.

\bibitem{mariere2003}
{\sc B.~Mariere, Z.-Q. Luo, and T.~N. Davidson}, {\em Blind constant modulus equalization via convex optimization}, IEEE Transactions on Signal Processing, 51 (2003), pp.~805--818.

\bibitem{marshall2008positive}
{\sc M.~Marshall}, {\em Positive polynomials and sums of squares}, no.~146, American Mathematical Soc., 2008.

\bibitem{nie2013certifying}
{\sc J.~Nie}, {\em Certifying convergence of lasserre’s hierarchy via flat truncation}, Mathematical Programming, 142 (2013), pp.~485--510.

\bibitem{nie2014optimality}
{\sc J.~Nie}, {\em Optimality conditions and finite convergence of {L}asserre's hierarchy}, Mathematical Programming, 146 (2014), pp.~97--121.

\bibitem{nie2023moment}
{\sc J.~Nie}, {\em Moment and polynomial optimization}, SIAM, 2023.

\bibitem{sinjorgo2024cuts}
{\sc L.~Sinjorgo, R.~Sotirov, and M.~F. Anjos}, {\em Cuts and semidefinite liftings for the complex cut polytope}, Mathematical Programming,  (2024), pp.~1--50.

\bibitem{smale1981fundamental}
{\sc S.~Smale}, {\em The fundamental theorem of algebra and complexity theory}, Bulletin (New Series) of the American Mathematical Society, 4 (1981), pp.~1--36.

\bibitem{smale1998mathematical}
{\sc S.~Smale}, {\em Mathematical problems for the next century}, The mathematical intelligencer, 20 (1998), pp.~7--15.

\bibitem{waldspurger2015phase}
{\sc I.~Waldspurger, A.~d’Aspremont, and S.~Mallat}, {\em Phase recovery, maxcut and complex semidefinite programming}, Mathematical Programming, 149 (2015), pp.~47--81.

\bibitem{wang2026more}
{\sc J.~Wang}, {\em A more efficient reformulation of complex sdp as real sdp: J. wang}, Computational Optimization and Applications,  (2026), pp.~1--17.

\bibitem{wang2022exploiting}
{\sc J.~Wang and V.~Magron}, {\em Exploiting sparsity in complex polynomial optimization}, Journal of Optimization Theory and Applications, 192 (2022), pp.~335--359.

\bibitem{wang2025real}
{\sc J.~Wang and V.~Magron}, {\em A real moment-hsos hierarchy for complex polynomial optimization with real coefficients}, Computational Optimization and Applications, 90 (2025), pp.~53--75.

\bibitem{tssos1}
{\sc J.~Wang, V.~Magron, and J.-B. Lasserre}, {\em {TSSOS}: A moment-{SOS} hierarchy that exploits term sparsity}, SIAM Journal on Optimization, 31 (2021), pp.~30--58.

\bibitem{tssos3}
{\sc J.~Wang, V.~Magron, J.-B. Lasserre, and N.~H.~A. Mai}, {\em {{CS-TSSOS}: Correlative and term sparsity for large-scale polynomial optimization}}, arXiv:2005.02828,  (2020).

\bibitem{zhang2006complex}
{\sc S.~Zhang and Y.~Huang}, {\em Complex quadratic optimization and semidefinite programming}, SIAM Journal on Optimization, 16 (2006), pp.~871--890.

\end{thebibliography}
\end{document}